\documentclass[10pt]{article}

\usepackage[margin=1in]{geometry}

\usepackage{amsmath,amssymb,amsthm}
\usepackage{mathtools}
\usepackage{geometry}
\usepackage{bm} 
\usepackage{xcolor} 
\usepackage{mathrsfs} 
\usepackage{comment} 
\usepackage[hidelinks]{hyperref} 
\usepackage{cleveref}
\usepackage{tikz}
\usepackage{pgfplots}
\usepackage{subcaption}
\pgfplotsset{compat=1.18}
\usetikzlibrary{decorations.pathreplacing}

\crefname{equation}{equation}{equations}
\Crefname{equation}{Equation}{Equations}

\newcommand{\R}{\mathbb{R}}
\newcommand{\Rd}[1]{\mathbb{R}^{#1}}

\newcommand{\vzero}{\bm{0}}
\newcommand{\va}{\bm{a}}
\newcommand{\vb}{\bm{b}}
\newcommand{\vf}{\bm{f}}
\newcommand{\vF}{\bm{F}}
\newcommand{\vk}{\bm{k}}

\newcommand{\vu}{\bm{u}}
\newcommand{\vv}{\bm{v}}
\newcommand{\vx}{\bm{x}}
\newcommand{\vxp}{\bm{x}'}
\newcommand{\vy}{\bm{y}}
\newcommand{\vpot}{\vv}
\newcommand{\vpoth}{\vpot_{H}}
\newcommand{\vpotl}{\vpot_{L}}

\newcommand{\vga}{\bm{\gamma}}
\newcommand{\veta}{\bm{\eta}}
\newcommand{\vmu}{\bm{\mu}}
\newcommand{\vnu}{\bm{\nu}}
\newcommand{\vsigma}{\bm{\sigma}}
\newcommand{\vxi}{\bm{\xi}}
\DeclarePairedDelimiter{\abs}{\lvert}{\rvert}

\newcommand{\fhat}[1]{\widehat{#1}}

\newcommand{\fg}{\fhat{g}}
\newcommand{\fF}{\fhat{F}}

\newcommand{\fv}{\fhat{v}}
\newcommand{\fu}{\fhat{u}}
\newcommand{\ful}{\fhat{u}_{L}}
\newcommand{\fuh}{\fhat{u}_{H}}

\newcommand{\fM}{\fhat{M}}

\newcommand{\calD}{\mathcal{D}}
\newcommand{\calF}{\mathcal{F}}

\newcommand{\calS}{\mathcal{S}}
\newcommand{\calV}{\mathcal{V}}

\newcommand{\pv}{\operatorname{PV}}
\newcommand{\distop}[2]{\left\langle #1,#2\right\rangle}
\newcommand{\real}{\operatorname{Re}}

\newcommand{\dom}{\Omega}
\newcommand{\partdom}{\partial\Omega}
\newcommand{\bdry}{\Gamma}

\newcommand{\calSL}{\mathcal{S}_L}

\newcommand{\calDL}{\mathcal{D}_L}

\newcommand{\calVL}{\mathcal{V}_L}

\newcommand{\ul}{u_{L}}
\newcommand{\uh}{u_{H}}

\newcommand{\slp}{\calS}     \newcommand{\slpl}{\calSL}
\newcommand{\dlp}{\calD}    \newcommand{\dlpl}{\calDL}
\newcommand{\vp}{\calV}     \newcommand{\vpl}{\calVL} \newcommand{\vpldd}{\calVL^{\delta}} \newcommand{\vplpv}{\calVL^{\pv}}

\newcommand{\ind}{\chi_{\dom}}
\newcommand{\veps}{\varepsilon}
\newcommand{\vphi}{\varphi}
\newcommand{\sgn}{\rho}
\newcommand{\ordo}[2]{\mathcal{O}\!\left(#1^{#2}\right)}
\newcommand{\curve}{\gamma}
\newcommand{\curvex}{\tilde{\curve}}
\newcommand{\schwartz}{\mathscr{S}}
\newcommand{\pschwartz}{\mathscr{S}'}

\DeclareMathOperator{\erf}{erf}
\DeclareMathOperator{\erfc}{erfc}

\title{Fourier-based potential theory without an explicit Green's function}
\author{Fredrik Fryklund%
  \thanks{Department of Mathematics, KTH Royal Institute of Technology, Stockholm, Sweden,
  ({\tt ffry@kth.se}).}}
\date{}

\begin{document}

\maketitle

\begin{abstract}
Integral equation methods provide an effective framework for solving partial differential equations, but their applicability typically relies on the availability of explicit free-space Green’s functions. For coupled systems arising in multiphysics applications, such Green’s functions are generally not available, limiting the scope of classical potential theory-based approaches. In this work, we introduce a formulation of potential theory that avoids explicit use of Green’s functions entirely, relying instead on the Fourier symbol of the governing operator. The central idea is a parabolic regularization of the symbol, which yields a decomposition of the solution into a smooth, nonlocal component and a spatially localized residual. For the localized component, we derive explicit asymptotic expansions for volume, single-layer, and double-layer potentials in powers of a length scale parameter $\veps$. The coefficients are expressed in terms of local geometric quantities and derivatives of the source data. The derivation is carried out entirely in the Fourier domain and applies to the Poisson equation in two and three dimensions, as well as to a class of coupled strongly elliptic systems.

\end{abstract}

\section{Introduction}\label{sec:introduction}

Two classical tools for the solution of constant-coefficient partial
differential equations (PDEs) are Fourier analysis and potential theory.
The former provides an analytic solution to inhomogeneous
PDEs in free-space or separable geometries, while the
latter provides powerful solution methods for homogeneous
PDEs in complicated domains.  Upon discretization,
both Fourier and integral equation methods lead to dense matrices, and 
the last decades have seen substantial progress on
fast algorithms that permit optimal (near linear scaling) solution methods, 
including the fast Fourier
transform (FFT), hierarchical methods such as tree codes,
the fast multipole method (FMM) and $\mathcal{H}$-matrices, the non-uniform FFT,
and linear algebraic methods that compress blocks of the system matrix that
corresponding to well-separated degrees of freedom \cite{hackbusch1999,LINDBO20118744,fds7,carrier1988sisc,greengard2004nufft}.    

Historically, however, Fourier methods have been viewed as unsuitable for 
handling highly nonuniform source distributions or complicated domains.
Potential theory, on the other hand, requires knowledge of the Green's function
corresponding to the equation of interest in order to define single and double
layer potentials, whose source distributions are restricted to the domain boundary.
While the Green's function is known in closed form for many of the equations of 
classical physics, that is not the case for coupled systems of the type that arise, 
for example, in multiphysics applications.

Here, we show that we can develop integral equation methods
in general geometries without access to the physical space Green's function,
using only the Fourier transform of the PDE (its {\it symbol}). That is, we show that asymptotic expansions can be developed for volume integrals with smooth or discontinuous data, as well
as layer potentials, {\it directly in the Fourier domain}.
This is the principal result of the present paper, which can be combined with other fast algorithms to yield solvers that are high-order accurate and fully adaptive, and are also robust to the quality of the underlying volumetric mesh, while achieving near-optimal computational scaling.

Before describing the scheme in detail, let us illustrate the basic idea
with a simple example. For this, we use the convention that, for a function 
$F(\vx)$, its Fourier transform is defined as
\begin{equation}\label{eq:fourierdef}
 \fF(\vxi) = \calF[F](\vxi) := \int_{\Rd{d}}F(\vx)e^{-i\vxi\cdot \vx}\,d\vx.
\end{equation}
We assume that 
the Fourier inversion formula 
\begin{equation}\label{eq:fourierinv}
 F(\vx) = \calF^{-1}[\fF](\vx) := 
\frac{1}{(2\pi)^d}\int_{\Rd{d}}\fF(\vxi)e^{i\vxi\cdot \vx}\,d\vxi
\end{equation}
also holds, which requires some conditions on $F$, $\fF$ that we postpone
until \cref{sec:regstrat}. Suppose now that we seek to solve the Poisson equation
\begin{equation}\label{eq:simplesys}
    -\Delta u(\vx) = g(\vx) 
\end{equation}
for $\vx \in \mathbb{R}^2$,
where $\vx = (x_1,x_2)$ and $\Delta$ is the Laplacian,
with suitable conditions at infinity so that the solution can be expressed
via the Fourier transform.
Then, after transforming both sides of \eqref{eq:simplesys},
we have
\begin{equation}\label{eq:simplesol}
u(\vx) = 
\frac{1}{(2\pi)^2}\int_{\Rd{2}}
\frac{\fg(\vxi)}{\abs{\vxi}^2} e^{i\vxi\cdot \vx}\,d\vxi.
\end{equation}
If $g(\vx)$ were smooth, compactly supported, and well-resolved by a uniform 
mesh, then \eqref{eq:simplesol} can be evaluated efficiently using the FFT
with care taken to handle the singularity at $\vxi=\vzero$ \cite{vico}.
When $g(\vx)$ is not smooth, however, the Fourier transform is slowly 
decaying and \eqref{eq:simplesol} leads to a less effective scheme.
Instead, we can write
\begin{equation}\label{eq:simplesplit}
u(\vx) = 
\frac{1}{(2\pi)^2}\int_{\Rd{2}}
\fg(\vxi)
\frac{e^{-\abs{\vxi}^2 \veps}}{\abs{\vxi}^2} 
e^{i\vxi\cdot \vx}\,d\vxi +
\frac{1}{(2\pi)^2}\int_{\Rd{2}}
\fg(\vxi)
\frac{1 - e^{-\abs{\vxi}^2 \veps}}{\abs{\vxi}^2} 
e^{i\vxi\cdot \vx}\,d\vxi,
\end{equation}
with the first part, denoted $u_H$, is a smoothed (mollified) solution,
and the second part $\ul$ can be thought of as the residual correction.
The smoothness of $u_H$ is enforced
by the factor $e^{-\abs{\vxi}^2 \veps}$, independent of any
assumptions about $g(\vx)$, but it is spatially nonlocal.
For nonuniform distributions requiring adaptive discretization,
fast algorithms for evaluating $u_H$ are now readily available,
essentially by replacing the FFT with the DMK method \cite{dmk2024}.
While the
residual correction $\ul$ is non-smooth, it is rapidly decaying
 in physical space.
This localization of $\ul$ 
follows from the Paley-Wiener theorem \cite{hormander1976linear} and the fact that
\[
\hat{G}_L(\vxi) = \frac{1 - e^{-\abs{\vxi}^2 \veps}}{\abs{\vxi}^2} 
\]
is infinitely differentiable.
The residual kernel in physical space is defined as
the inverse Fourier transform of $\hat{G}_L(\vxi)$.
In previous work \cite{lghtwght_ps2026}, we showed that
the convolution of the residual kernel
with smooth, discontinuous or distributional data (layer potential densities)
can be carried out to high precision using asymptotic methods.
It turns out, as we shall see below,
that these asymptotic approximations can be derived directly from 
the Fourier transform for a large class of PDEs, permitting a significant 
generalization of the 
results of \cite{lghtwght_ps2026} which were based on the explicit form of the 
residual kernel in physical space. The simplest such example is the case when
$g(\vx)$ is globally smooth. Then, 
\begin{align}\label{eq:simplesasymp}
\ul(\vx) &= 
\frac{1}{(2\pi)^2}\int_{\Rd{2}}
\fg(\vxi)
\left[ \frac{1 - e^{-\abs{\vxi}^2 \veps}}{\abs{\vxi}^2} \right]
e^{i\vxi\cdot \vx}
\,d\vxi  \notag \\
&=
\frac{1}{(2\pi)^2}\int_{\Rd{2}}
{\fg(\vxi)}
\left[ \veps - (1/2)\veps^2 \abs{\vxi}^2 + O(\veps^3) \right]
e^{i\vxi\cdot \vx}
\,d\vxi,
\end{align}
after Taylor expansion of $e^{-\abs{\vxi}^2 \veps}$.
Taking the inverse Fourier transform leads directly to the formula
\begin{equation}\label{eq:simplesasymp2}
\ul(\vx) = 
\veps g(\vx) +
\frac{\veps^2}{2} \Delta g(\vx) + O(\veps^3).
\end{equation}

Suppose now that $g(\vx)$ is discontinuous, with 
$g(\vx) = e^{-(x_1^2 +x_2^2)}$ for $x_2>0$ and
$g(\vx) = 0$ for $x_2\leq 0$. There is no significant difference 
in the treatment of the regularized part. For the local part,
$\ul$, however, we need to exploit the spatially localized kernel $S_L$ as a function of $\veps$. 
Noting that $S_L= \int_0^\veps e^{-\abs{\vxi}^2 t} \,dt$, 
and using the Fourier transform pair 
$\calF[e^{-t\abs{\xi}^2}] \leftrightarrow \sqrt{\tfrac{\pi}{t}}e^{-\abs{x}^2/(4t)}$,
we may write \eqref{eq:simplesasymp} as
\begin{align}
  \ul(\vx) = \frac{1}{4\pi}\int_0^{\veps}\int_{\R}\int_{0}^{\infty}
             g(\vxp)
             e^{-\abs{\vxp-\vx}^2/(4t)}\,dx_2'\,dx_1'\,dt  
\end{align}
from Fubini's theorem. We now assume the target point $\vx$ is located at $\vx=(0,r)$, and Taylor expand $g(\vxp)$ to second order around $(0,r)$, followed by integration in the tangential variable $x_1$ over the real line, giving
\begin{equation}
  \ul(\vx) = \frac{1}{2\sqrt{\pi}}\int_0^{\veps}t^{-1/2}\int_{-r}^{\infty}  e^{-(x'_2)^2/(4t)}\left[g(0,r) + x'_2g_{x'_2}(0,r) + \ordo{(x'_2)}{2}\right]\,dx'_2\,dt.
\end{equation}
Evaluating the resulting Gaussian integrals explicitly yields
\begin{align}\label{eq:halfspaceasymp}
  \ul(\vx) &= 
             \veps A_0(c)\,g(0,r)
             +
             \veps^{3/2}A_1(c)\,g_{x'_2}(0,r)
             +
             \mathcal O(\veps^2),
\end{align}
with $c=\frac{r}{\sqrt{\veps}}$, where the coefficients are
\begin{align}
  A_0(c)
  =
  \frac12\left(
  1+\operatorname{erf}\!\left(\frac c2\right)
  +\frac{c}{\sqrt{\pi}}e^{-c^2/4}
  -\frac{c^2}{2}\operatorname{erfc}\!\left(\frac c2\right)
  \right)\qquad
A_1(c)
=
\frac{c^3}{6}\operatorname{erfc}\!\left(\frac c2\right)
+
\frac{2-c^2}{3\sqrt{\pi}}e^{-c^2/4}
\end{align}
and $\erf(x)$ denotes the usual error function.
We note that \eqref{eq:halfspaceasymp} is reduced to \eqref{eq:simplesasymp2} when the distance r to the $x_2$-axis goes to infinity. The derivations become more involved for 
volume integrals over complicated domain boundaries and layer potentials, and will require the theory of tempered distributions, but the essential idea is the same.

The reader may have noted that, in \eqref{eq:simplesasymp2}, the asymptotic 
expansion involves only integer powers of $\veps$, while
in \eqref{eq:halfspaceasymp}, half-order powers appear. This follows from a careful
derivation of \eqref{eq:halfspaceasymp} in a manner consistent with the parabolic 
regularization strategy of \eqref{sec:regstrat}.

This paper is organized as follows. The main result is illustrated
through the canonical example of the free-space Poisson equation in two 
dimensions. In \Cref{sec:pottheory}, the necessary elements of potential theory are reviewed, 
followed by \Cref{sec:regstrat} where
the parabolic regularization strategy is developed. In \Cref{sec:mainresult} the asymptotic expansions of the local parts are derived, and \Cref{sec:poissonr3} and \Cref{sec:coupsys} extend the construction to the Poisson equation in three dimensions and to coupled elliptic system, respectively.

\section{Potential theory}\label{sec:pottheory}
 
There are a number of advantages of recasting the
PDEs of mathematical 
physics as integral equations. The governing Green's function 
expresses outgoing radiation conditions exactly so that 
there is no need for artificial boundary conditions on
a finite computational domain.
In the vector case (say, for the Stokes or Maxwell equations)
it can also enforce divergence-free constraints on velocity 
or electromagnetic fields, avoiding the need to impose such 
constraints algebraically. Integral equation methods
also reduce the conditioning of
the problem, by restricting unknowns to the domain
boundary, even for inhomogeneous equations with volumetric 
sources. Finally, when carefully designed, they lead to 
well-conditioned linear systems of equations and more 
robust {\it a priori} and {\it a posteriori} error
estimates. 

Given a Green's function, 
a rather general representation of a potential function 
$u(\vx)$ is
\begin{equation}\label{eq:pot_basic_split}
  u(\vx) = \int_\dom G(\vx-\vxp)f(\vxp)\,d\vxp
          + \int_{\partdom} G(\vx-\vxp)\sigma(\vxp)\,ds_{\vxp}
          + \int_{\partdom} \frac{\partial G(\vx-\vxp)}{\partial\vnu_{\vxp}}
            \mu(\vxp)\,ds_{\vxp}.
\end{equation}
(In the scalar case, recall that any $C^2$ function can be expressed
in this form by Green's third identity.)
The first integral operator above is referred to as a volume
potential, the second as a single
layer potential and third as a double layer potential. 
In order for such operators to be of practical interest,
two main  computational challenges need to be addressed.
First, $G$ can decay slowly, so that 
discretization leads to dense matrices.
Second, $G$ is (weakly) singular, requiring sophisticated 
quadrature methods for accurate evaluation. 

Historically, these challenges have been
addressed in one of two ways.
One can decompose $u(\vx) = u_F(\vx) + u_N(\vx)$, where $u_F$ 
corresponds to ``well-separated" far-field interactions
and $u_N$ corresponds to 
the near field, ranging over the sources or source distribution
in some small neighborhood of the target $\vx$.
The far field interactions are smooth so that $u_F$ can be 
computed using standard quadrature rules for smooth functions.
The near field requires more sophisticated quadrature rules 
that typically integrate singular kernels multiplied by 
smooth densities over discrete volumetric or surface patches.
An alternative approach is to write
$G = G_H + G_L$,
where $G_H$ is a mollified but still nonlocal version of $G$ and 
$G_L$ is a localized correction.
This is distinct from the near/far decomposition
$u(\vx) = u_F(\vx) + u_N(\vx)$ above, where
$G$ is unaltered.
Such a splitting is the basis for Ewald summation 
in $\mathbb{R}^3$:
\[ \frac{1}{r} = \frac{\erf(r/\veps)}{r} + 
\frac{\erfc(r/\veps)}{r},
\]
where $\erf$ and $\erfc$ are the classical error and 
complementary error functions, respectively.

Fast algorithms reduce the cost of 
evaluating~\cref{eq:pot_basic_split} with $N$ discretization points (and targets)
from $O(N^2)$ to $O(N)$ or $O(N\log N)$. Adaptive approaches that rely on near/far splitting
are generally based on hierarchical compression of 
well-separated interactions such as the fast multipole method  
(see, for example, \cite{carrier1988sisc, hackbusch1999, kenho1, fds7}).
For quasiuniform discretizations, the pre-corrected FFT,
spectral Ewald summation, and related methods 
project data onto a Cartesian grid to exploit the 
speed of the FFT~\cite{LINDBO20118744} for computing convolutions with the
regularized part ($\frac{\erf(r/\veps)}{r}$), while the residual kernel
$(\frac{\erfc(r/\veps)}{r})$ is handled in physical space.
For problems with multiple spatial scales,
the dual-space multilevel kernel-splitting (DMK) method
of \cite{dmk2024} permits linear-scaling Fourier-based convolution 
in a multiscale framework by
combining ideas from tree-based methods, multilevel summation
and Ewald summation. We do not seek to review the literature here and
refer the reader to \cite{dmk2024} for a historical overview.

Once the far field or mollified interactions have been accounted
for, using either hierarchical compression or DMK,
it remains to compute a convolution with the residual kernel.
For volume and layer potentials, this still leaves open the 
problem of quadrature with singular or weakly singular
integrands. This is an active area
of research and
a large variety of methods that are compatible with 
fast algorithms have been developed for this purpose, including
singularity subtraction  \cite{helsingojala2008}, generalized Gaussian quadrature \cite{ggq2010}, Quadrature By Expansion
(QBX) \cite{qbx2012} , etc.
We refer the reader to \cite{dmk2024,lghtwght_ps2026} and the references therein. 

Using the Ewald-like kernel splitting approach,
asymptotic methods have also been introduced
for integrating the residual kernel.
For the Laplace equation, the kernel 
$G_L = \frac{\erfc(r/\veps)}{r}$ is closely connected
to the heat equation and this correspondence was used in 
~\cite{greengardstrain1990,strain1992,JingLi,wang2025}
to develop expansions in powers of $\veps$.
These were carried out to higher order and to volume
potentials in \cite{lghtwght_ps2026}.
Related asymptotic analysis has been carried out by Beale,
Tlupova et al. (see, for example, \cite{beale2024}).
In~\cite{beale2025}, an interesting alternative 
was proposed, where the mollified kernel itself
is altered using asymptotic analysis so that, for a layer potential, 
integration against a smooth density leads to high order 
accuracy with respect to a length scale parameter $\veps$.
Somewhat different asymptotic methods were developed in~\cite{khatri2020} 
for nearly singular integrals to accelerate
the evaluation of field close to a given surface.

In the present work, we describe a parabolic regularization strategy,
motivated by the work in \cite{lghtwght_ps2026}, that
can be carried out for any strongly elliptic system 
in the Fourier transform domain once the symbol is known.
We show how to treat the single layer, double layer, and volume potentials based on the residual kernel.
This requires some care and the use of
tempered distributions, leading to asymptotic expansions for the local
parts in powers of $\sqrt{\veps}$, with coefficients given by geometric data and source derivatives at the target location.
We demonstrate the derivation and performance of the method
for the Poisson equation in two and three dimensions, and for a 
coupled strongly elliptic system. 

\section{The parabolic regularization strategy}\label{sec:regstrat}
In this section, we develop the parabolic regularization strategy for the free-space Poisson equation in two dimensions with compactly supported distributional source terms. This model problem is simple enough to make the regularization procedure transparent, while still containing the key features needed for the general setting. 
Our starting point is a representation of the elliptic solution as the steady-state
solution to an auxiliary heat equation, which we then split into its short time 
and long time components based on a parameter
$\veps>0$. The local part $\ul$ is the short time component, which captures the 
elliptic Green's function 
singularity and is exponentially localized. The smooth, regularized part $\uh$
corresponds to the long time component or ``history".
As noted in the introduction, the critical point is that the decomposition is 
constructed from the Fourier symbol of the operator rather than from the
physical-space Green's function itself. This will permit us to employ
the same strategy for more general strongly elliptic systems.

\subsection{Setup}
Consider the free-space Poisson equation in the plane,
\begin{equation}\label{eq:ps_2d}
   -\Delta u = F \quad \text{in }\Rd{2}.
\end{equation}
where $F$ has compact support. Also, $F$ is a tempered distribution, denoted $F\in\pschwartz(\Rd{2})$ which is the space of continuous linear functionals on the Schwartz space $\schwartz(\Rd{2})$, consisting of smooth functions on $\Rd{2}$ whose derivatives decay faster than any polynomial. The Poisson equation is understood in the distributional sense:
\begin{equation}\label{eq:ps_2d_dist}
\distop{-\Delta u}{\vphi}=\distop{F}{\vphi}
\quad \text{for all }\vphi\in \schwartz(\Rd{2})
\end{equation}
meaning $\distop{T}{\vphi}$ denotes the action of a distribution $T$ on a test function $\vphi$. As an additional criterion, we impose $\distop{F}{1} = 0$.

The reason for working at this level of generality is that the regularization strategy must apply not only to smooth volumetric data, but also to singular sources such as layer potentials. For this reason we formulate the problem in the language of distributions from the outset and work at the level of distributions throughout this section.

Recall that if $T\in\pschwartz(\Rd{2})$ is a tempered distribution, its Fourier transform $\fhat{T}\in\pschwartz(\Rd{2})$ is defined by duality,
\begin{equation}
  \langle\fhat{T},\vphi\rangle
  = \distop{T}{\fhat{\vphi}},
  \quad \text{for all }\vphi\in\schwartz(\Rd{2}),
\end{equation}
where $\fhat{\vphi}$ denotes the classical Fourier transform of the Schwartz function $\vphi$. Since the source term $F\in\pschwartz(\R^2)$ in \cref{eq:ps_2d} has compact support, it can be paired with any smooth function, so for every $\vxi\in\Rd{2}$ the pairing $\distop{F}{e^{-i\vxi\cdot\vx}}$ is well defined. In particular, the Fourier transform of $F$ may be written pointwise as
\begin{equation}
\fF(\vxi)=\distop{F}{e^{-i\vxi\cdot\vx}}, \quad \vxi\in\Rd{2}.
\end{equation}

In two dimensions, solutions of \cref{eq:ps_2d} are determined only up to an additive constant and may have far-field behavior of the form
\begin{equation}
u(\vx)\sim -\frac{\fF(\vzero)}{2\pi}\log\abs{\vx},
\quad\text{as }\abs{\vx}\to\infty,
\end{equation}
where $\fF(\vzero) = \distop{F}{1}$. Thus we remove the logarithmic growth at infinity by imposing $\fF(\vzero) = 0$ and we obtain uniqueness by fixing the additive constant to zero.

\subsection{The auxiliary parabolic problem}
The splitting of $u$ into a localized part $\ul$ and a non-localized part $\uh$ is obtained by using a parabolic regularization, parameterized by $\veps >0 $. Consider the auxiliary evolution problem
\begin{equation}\label{eq:heateq}
\partial_t v - \Delta v = F, \quad \text{in }(0,\infty)\times\Rd{2},
\end{equation}
for a tempered distribution $v(t)\in\pschwartz(\Rd{2})$ continuous in $t$ with zero initial data $v(0)=0$ and with $F$ from \cref{eq:ps_2d}. We recognize \eqref{eq:heateq} as the heat equation and will show that its steady-state limit in $t$ is the solution $u$ to \cref{eq:ps_2d}.

Applying the Fourier transform to \cref{eq:heateq} gives, for each $\vxi\in\Rd{2}$, an ordinary differential equation (ODE)
\begin{equation}\label{eq:heateq_fourier}
  \partial_t \fv (t,\vxi) + \abs{\vxi}^2 \fv (t,\vxi) = \fF(\vxi), \quad \fv(0,\vxi) = 0.
\end{equation}
Since $\fF$ is smooth, \eqref{eq:heateq_fourier} holds in a pointwise sense for every fixed $\vxi\in\Rd{2}$. For $\vxi \neq \bm{0}$, the unique solution with zero initial data is
\begin{equation}\label{eq:heat_ode_sol}
 \fv(t,\vxi) = \fhat{F}(\vxi)\int_0^t e^{-\tau\abs{\vxi}^2}\,d\tau = \frac{1-e^{-t\abs{\vxi}^2}}{\abs{\vxi}^2}\fF(\vxi),\quad \abs{\vxi} \neq \bm{0}.
\end{equation}
At $\vxi=\vzero$, \cref{eq:heateq_fourier} reduces to $\partial_t\fv(t,\vzero)=\fF(\vzero)$ with $\fv(0,\vzero)=0$, whose solution is $\fv(t,\vzero)=\fF(\vzero)\,t$. By the normalization $\fF(\vzero)=\distop{F}{1}=0$, this gives $\fv(t,\vzero)=0$ for all $t\ge 0$.
\subsection{Steady-state limit}\label{ss:steadystatelim}
The limit of the function $v(t)$ in $\pschwartz(\R^2)$ as $t \to \infty$ is the solution $u$ of \cref{eq:ps_2d}. To see this, let $\vxi\neq \bm{0}$, then \eqref{eq:heat_ode_sol} gives
\begin{equation}
\fv(t,\vxi) = \frac{1-e^{-t\abs{\vxi}^2}}{\abs{\vxi}^2}\fF(\vxi)\to  \frac{\fF(\vxi)}{\abs{\vxi}^2} \quad \text{ as } t \to \infty,
\end{equation}
and at $\vxi = \vzero$ we have $\fhat{v}(t,\vzero) = 0$. We identify
\begin{equation}
  \fu(\vxi) := \frac{\fF(\vxi)}{\abs{\vxi}^2},\quad \text { for  }\vxi\neq\bm{0},\quad\text{ and  } \fu(\vzero)=0,
\end{equation}
where $\fu$ is defined pointwise in all of $\Rd{2}$. To express the difference between $\fhat{v}(t)$ and the limit $\fu$, and to show $v(t)\to u$ in $\pschwartz(\Rd{2})$, it is convenient to introduce the Gaussian multiplier
\begin{equation}\label{eq:m_t}
m_t(\vxi):=e^{-t\abs{\vxi}^2},
\end{equation}
which belongs to $\schwartz(\Rd{2})$ for each fixed $t>0$. From \eqref{eq:heat_ode_sol} and the definition of $\fu$, for $\vxi\neq\vzero$,
\begin{equation}\label{eq:vhat_minus_uhat}
  \fv(t,\vxi) - \fu(\vxi)
  = \frac{1-e^{-t\abs{\vxi}^2}}{\abs{\vxi}^2}\fF(\vxi)
    - \frac{\fF(\vxi)}{\abs{\vxi}^2}
  = -e^{-t\abs{\vxi}^2}\,\fu(\vxi)
  = -m_t(\vxi)\,\fu(\vxi),
\end{equation}
while both sides vanish at $\vxi=\vzero$, so \eqref{eq:vhat_minus_uhat} holds on all of $\Rd{2}$. Since $\fF\in C^\infty(\Rd{2})$ with $\fF(\vzero)=0$, the function $\fu=\fF/\abs{\vxi}^2$ has at most a $\abs{\vxi}^{-1}$ singularity at
the origin and is locally integrable. Together with the polynomial growth of $\fF$ at infinity, this gives $\fu\,\vphi\in L^1(\Rd{2})$ for every $\vphi\in\schwartz(\Rd{2})$, so $\fu$ acts on such $\vphi$ by integration.

Fix $\vphi\in\schwartz(\Rd{2})$. By the duality of the Fourier transform and
\eqref{eq:vhat_minus_uhat},
\begin{equation}
  \distop{v(t)-u}{\fhat{\vphi}}
  = \distop{\fv(t)-\fu}{\vphi}
  = -\distop{m_t\fu}{\vphi}
  = -\distop{\fu}{m_t\vphi}
  = -\int_{\Rd{2}} \fu(\vxi)\,m_t(\vxi)\,\vphi(\vxi)\,d\vxi,
\end{equation}
where we used that the definition of multiplication of a tempered distribution $T$ by a smooth function $m_t$ is $\distop{m_tT}{\vphi} = \distop{T}{m_t\vphi}$ for all $\vphi\in\schwartz(\Rd{2})$.

The integrand is bounded by $\abs{\fu\,\vphi}\in L^1(\Rd{2})$ uniformly in $t$ and tends to $0$ pointwise as $t\to\infty$, so $\distop{v(t)-u}{\fhat{\vphi}}\to 0$ by th dominated convergence theorem, and hence $v(t)\to u$ in $\pschwartz(\Rd{2})$ as $t\to\infty$. The steady-state limit of \cref{eq:heateq} is therefore the solution of \cref{eq:ps_2d}.

\subsection{The local and history parts}
From \cref{ss:steadystatelim}, we have that taking $t\to \infty$ in \cref{eq:heat_ode_sol} gives
\begin{equation}
  \fu (\vxi) = \fF(\vxi)\int_0^{\infty}e^{-\tau\abs{\vxi}^2}\,d\tau,\qquad \vxi\neq \vzero,
\end{equation}
with $\fu(\vzero)=0$. For $\veps>0$, we split this integral at $\tau = \veps$ and define the local part 
$\ul$ of $u$ in Fourier space by
\begin{equation}\label{eq:pot_ps2d_floc}
\ful(\vxi) := \fF(\vxi)\int_0^\veps e^{-\tau\abs{\vxi}^2}\,d\tau.
\end{equation}
The history part $\uh$ of $u$ is defined in Fourier space by
\begin{equation}\label{eq:pot_ps2d_fhis}
\fuh(\vxi) := \fF(\vxi)\int_\veps^\infty e^{-\tau\abs{\vxi}^2}\,d\tau,\quad \vxi\neq\vzero,
\end{equation}
with $\fuh(\vzero) := 0$, ensuring $\fu=\ful+\fuh$ holds on all of $\Rd{2}$. 

Taking the inverse Fourier transform of \cref{eq:pot_ps2d_floc,eq:pot_ps2d_fhis}, we write $\ul = \calF^{-1}[\ful]$ and $\uh = \calF^{-1}[\fuh]$. For the local part, our goal is now to rewrite the representation so that the $\tau$–integral appears as the outermost integral, reflecting the underlying parabolic evolution and providing the form required for the asymptotic expansions in \cref{sec:mainresult}. 

The source density $F$ has compact support, thus the Paley--Wiener theorem \cite{hormander1976linear} implies that $\fF\in C^\infty(\Rd{2})$ and that all its derivatives have at most polynomial growth. It follows that $\fF\,\fhat{\phi}\in L^1(\Rd{2})$ for any $\phi\in\schwartz(\Rd{2})$. Moreover, for each $\tau>0$ we have
\begin{equation}
\distop{m_\tau\fF}{\fhat{\phi}} = \int_{\Rd{2}} e^{-\tau\abs{\vxi}^2}\fF(\vxi)\fhat{\phi}(\vxi)\,d\vxi,
\end{equation}
with $m_\tau$ from \eqref{eq:m_t}. Since the integral is finite, independently of $\tau$, and the map $\tau\mapsto \distop{m_\tau\fF}{\fhat{\phi}}$ belongs to $L^1([0,\veps])$,  the distribution-valued integral
\begin{equation}
\distop{\int_0^\veps m_\tau\fF\,d\tau}{\fhat{\phi}}
= \int_0^\veps \distop{m_\tau\fF}{\fhat{\phi}}\,d\tau
\end{equation}
is well-defined. Using the duality of the Fourier transform, we obtain
\begin{equation}\label{eq:pot_ps2d_loc_deriv}
\distop{\calF^{-1}[\ful]}{\phi} = \int_0^\veps \distop{m_\tau\fF}{\calF^{-1}\phi}\,d\tau= \int_0^\veps \distop{\calF^{-1}[m_\tau\fF]}{\phi}\,d\tau.
\end{equation}
For each $\tau>0$, the function $m_\tau\fF$ is rapidly decaying together with all its derivatives, and belongs to $\schwartz(\Rd{2})$. It follows that $\calF^{-1}[m_\tau\fF]$ is given by a classical absolutely convergent integral for every $\vx\in\Rd{2}$, and depends continuously on $\tau>0$. Consequently, the outer $\tau$–integral is an ordinary integral for each fixed $\vx\in \R^2$, whenever the integrand is integrable on $(0,\veps]$. For the sources considered here, namely piecewise smooth volume densities in bounded domains and layer distributions supported on smooth curves, this integrability holds at every interior target, and we obtain the pointwise representation
$u[F](\vx)=\ul[F](\vx)+\uh[F](\vx)$ with local part
\begin{equation}\label{eq:uL_def}
      \ul[F](\vx) = \frac{1}{(2\pi)^{d}} \int_0^{\veps} \int_{\Rd{d}} e^{-t\abs{\vxi}^2}\distop{F}{e^{-i \vxi \cdot \vxp}} \, e^{i \vxi \cdot \vx}\, d\vxi \, dt,
\end{equation}
and history part
\begin{equation}\label{eq:uH_def}
    \uh[F](\vx) = \frac{1}{(2\pi)^{d}} \int_{\veps}^{\infty} \int_{\Rd{d}} e^{-t\abs{\vxi}^2}\distop{F}{e^{-i \vxi \cdot \vxp}}\, e^{i \vxi \cdot \vx}\, d\vxi \, dt,
\end{equation}
where we keep $d$ general to reuse these in \cref{sec:poissonr3}. This holds for compactly supported source distributions $F$ generated by smooth densities on smooth geometric supports. More precisely, the sources are piecewise smooth volume densities in bounded domains and layer distributions supported on smooth curves. In the present paper this includes volume sources, single layer sources, and double layer sources.

In terms of the heat equation \eqref{eq:heateq}, the local part corresponds to the short-time heat evolution over $[0,\veps]$. At time $\tau$, the heat kernel is proportional to $e^{-\abs{\vx}^2/4\tau}$, and thus concentrates near the origin as $\tau\to 0$. Integrating over $\tau\in(0,\veps]$ therefore produces a contribution that is exponentially localized in physical space on a length scale of order $\sqrt{\veps}$, and it is precisely this term that captures the singularity of the free-space Green's function $G$. The history part, by contrast, corresponds to the evolution over $[\veps,\infty)$. In this regime the heat kernel is no longer concentrated, so the resulting contribution remains spatially smooth and represents the long-time evolution of the solution away from the singularity.

\section{Asymptotics for the Poisson equation in $\R^2$}\label{sec:mainresult}

We first illustrate the regularization strategy from
\cref{sec:regstrat} in the setting of the Poisson equation in $\R^2$.
Our objective is to derive explicit asymptotic expansions of the local part $\ul$
for three source types: a single layer density $\sigma$ supported on a curve $\bdry$,
a double layer density $\mu$ supported on $\bdry$, and a volume source density $f$
supported on a compact domain $\dom$. More precisely, these sources define the corresponding single layer, double layer, and volume potentials:
\begin{align}
  \label{eq:single_layer}
  \slp[\sigma](\vx) &= -\frac{1}{2\pi} \int_\bdry \log|\vx-\vy| \,\sigma(\vy)\, ds(\vy),\\
  \label{eq:double_layer}
  \dlp[\mu](\vx) &= \frac{1}{2\pi} \int_\bdry
  \frac{(\vx - \vy)\cdot \vnu(\vy)}{|\vx - \vy|^2} \,\mu(\vy)\, ds(\vy),\\
  \label{eq:volume_layer}
  \vp[f](\vx) &= -\frac{1}{2\pi} \int_\dom \log|\vx-\vy| \,f(\vy)\, d\vy.
\end{align}
in terms of the (known) free-space Green's function
\begin{equation}
G(\vx) = -\frac{1}{2\pi}\log\abs{\vx}.
\end{equation}

Rather than evaluating these expressions directly, we compute $\slp$, $\dlp$, and $\vp$ through the decomposition of the solution into a local part \eqref{eq:uL_def} and a history part \eqref{eq:uH_def}, obtained from the Fourier representation with sources of the form $F=\sigma\delta_\bdry$, $F=-\nabla\cdot(\vnu\mu\delta_\bdry)$ or $F=f\ind$, where $\ind$ denotes the indicator of the domain $\dom$ and $\delta_\bdry$ the restriction to the curve $\bdry$. The analysis therefore reduces to understanding the local contribution $\ul$ associated with these sources.

We derive asymptotic formulas of the kind
\begin{equation}
  \ul[F](\vx) = \sum_{p = 0}^{P}\veps^{p/2} A_p(\vx) 
  + \ordo{\veps}{(P+1)/2},
\end{equation}
for the local parts, where $A_p$ are coefficients containing geometric data and derivatives of $F$. In applications, $\veps = \ordo{h}{2}$ where $h$ is the spatial resolution, thus including $P$ terms gives a $(P+1)$th order scheme under refinement.

We consider the history part formulated as
\begin{equation}
\uh[F](\vx) = \frac{1}{(2\pi)^{d}} \int_{\Rd{d}} \frac{e^{-\veps\abs{\vxi}^2}}{\abs{\vxi}^2} \fF(\vxi)\, e^{i \vxi \cdot \vx}\, d\vxi.
\end{equation}
This representation of $\uh$ is equivalent to the time-integrated form, but is written in a form suitable for numerical computations. Since it is smooth and nonlocal, it can be evaluated to high accuracy using either FFT-based methods when the source $F$ is resolved on a uniform mesh, or using the DMK algorithm when multiscale features are present. We do not pursue any analysis  of the history part further, and instead focus on the local part.

We treat the single layer potential in \cref{subsec:slppois}, the double layer potential in \cref{subsec:dlppois}, and the volume potential in \cref{subsec:vppois}. While the geometric setting differs in the volume case, where the source is supported in a region rather than on a curve, the resulting asymptotic structure is analogous. In all three cases, the formulas obtained agree with those derived using physical space techniques in \cite{lghtwght_ps2026}.

\subsection{Geometric setting}\label{subsec:geometricsetting}

We now describe the geometric setting used in the asymptotic derivations.
Let $\dom \subset \R^2$ be a domain with boundary $\bdry = \partial \dom$
of class $C^k$, where $k$ is sufficiently large to support the required
derivatives. We assume that $\bdry$ admits a local parameterization
with nonvanishing Jacobian and that it can also be
represented as the zero level set of a function $g$, so that
\begin{equation}
\bdry = \{\vx \in \R^2 : g(\vx) = 0\},
\end{equation}
where $g$ is the signed distance function.

More precisely, and to
carry out the asymptotic expansions, we introduce local coordinates
centered at the closest boundary point $\vb$ to an arbitrary target
point $\vx \in \R^2$. We translate $\vb$ to the origin and rotate the coordinates so that the $x_1$-axis is tangent
to $\bdry$ at $\vb$ and the $x_2$-axis is aligned with the inward normal.
We assume that the boundary is represented locally as the graph
\begin{equation}\label{eq:loc_param}
\vga(x_1) = (x_1, \curve(x_1)), \qquad x_1 \in I_b,
\end{equation}
for a small interval $I_b$. By assumption, it satisfies $\curve(0)=\vb$, $\curve'(0)=0$ and $\curve''(0)=\kappa_{b}$ which is the curvature at $\vb$.

In these coordinates, the target point is at $\vx = (0,r)$, where
$r = |\vx - \vb|$ is the distance to the boundary. The relative position can be written as
\begin{equation}\label{eq:sgn_def}
\vx - \vb = \sgn\, r\, \vnu, \qquad
\sgn = \operatorname{sign}\big((\vx-\vb)\cdot \vnu\big) \in \{-1,0,1\},
\end{equation}
where $\vnu = \vnu(\vb)$ is the outward unit normal at $\vb$, as shown in \cref{fig:locparam}. For all derivations we assume $\vx$ is an interior point.
This local coordinate system will be used to reduce the
integrals defining the local part $\ul$ to a form suitable for asymptotic
analysis.

\begin{figure}
\centering
\begin{tikzpicture}

  \draw[thick] plot[domain=-1.8:2.8, samples=100]
    (\x, {-0.2*\x*\x + 0.08*\x*\x*\x})
    node[right] {$\curve(x_1)$};

  \fill (0,0) circle (2pt) node[above right] {$\bm{b}$};

  \draw[->, thick] (0,0) -- (1,0) node[above=4pt, right=-4pt] {$x_1$};

  \draw[->, thick, dashed] (0,0) -- (0,-1) node[right] {$\vnu$};

  \def\r{1.8}
  \fill (0,\r) circle (2pt) node[right=3pt] {$\bm{x}$};

  \draw[decorate, decoration={brace, amplitude=5pt, raise=4pt}]
    (0,0) -- (0,\r) node[midway, left=10pt] {$r$};

  \end{tikzpicture}
  \caption{Local coordinate system at a boundary point $\bm{b}$. The boundary is represented as a graph $\curve(x_1)$ with tangent direction $x_1$ and outward unit normal $\vnu$. The target point $\bm{x}$ lies at distance $r$ along the inward normal direction $-\vnu$ from $\bm{b}$.}
  \label{fig:locparam}
\end{figure}

\subsection{The single layer potential $\slp$}\label{subsec:slppois}
We first consider the single layer potential $\slp$ from \eqref{eq:single_layer},
with density $\sigma$ supported on a closed curve $\bdry$ with the goal of
computing the local part $\slpl$,
defined by \eqref{eq:uL_def}, with $F(\vx) = \sigma(\vx)\delta(g(\vx))$:
\begin{equation}\label{eq:slpl}
  \slpl[\sigma](\vx) := \ul[\sigma \delta(g)](\vx)
  = \frac{1}{(2\pi)^2} \int_0^{\veps} \int_{\R^2} e^{-t\abs{\vxi}^2}
  \int_{\R^2} \sigma(\vxp)\delta(g(\vxp))\, e^{-i \vxi \cdot \vxp} \, d\vxp
  \, e^{i \vxi \cdot \vx}\, d\vxi \, dt.
\end{equation}
The resulting expansion will take the form
\begin{equation}\label{eq:slplgeneral}
  \slpl[\sigma](\vx)
  = \sum_{p = 0}^{P}\veps^{p/2} A_p(\vx) 
  + \ordo{\veps}{(P+1)/2},
\end{equation}
where the coefficients $A_p(\vx)$ depend on $\sigma$, the local geometry of
$\bdry$, and their derivatives evaluated at the boundary point $\vb$ closest
to $\vx$. 

\subsubsection{$\slpl$ in the local coordinate system}\label{subsubsec:slppois_locparam}
In the derivation below we assume $r>0$, so that the target point is separated from the source curve. Then, after pairing the Fourier expression with the Gaussian factor, the resulting $z$-integrand is continuous on $(0,\sqrt\veps]$.

For sufficiently small $\sqrt{\veps}$, we restrict our attention to a neighborhood $U$ of the closest boundary point $\vb$ where the local graph representation \eqref{eq:loc_param} is valid. Choose a neighborhood $B$ of $\vb$, and a smooth cutoff function $\Psi\in C_c^{\infty}(\Rd{2})$ supported in $U$ such that $\Psi=1$ in $B$ and $\bar{B}\subset U$. Then the contribution of $\slpl$ from $\bdry\setminus U$ where the local graph is not valid, is of $\ordo{e}{-d_B^2/(4\veps)}$, where $d_B$ is the shortest distance from $\vx$ to $\overline{\bdry\setminus B}$. To see this, consider $\slpl[\sigma](\vx)~=~\slpl[\Psi\sigma](\vx)+\slpl[(1~-~\Psi)\sigma](\vx)$ where
  \begin{equation}\label{eq:slp_err_local}
    \begin{split}
      \slpl[(1-\Psi)\sigma](\vx) &=\frac{1}{(2\pi)^2}\int_{0}^{\veps}\int_{\Rd{2}}e^{-\abs{\vk}^2 t}\int_{\Rd{2}}\sigma(\vxp)\delta(g(\vxp))(1-\Psi(\vxp))e^{-i\vk\cdot(\vxp-\vx)}\,d\vxp d\vk dt\\
                                 &=\int_{0}^{\veps}\int_{\Rd{2}}\sigma(\vxp)\delta(g(\vxp))(1-\Psi(\vxp))\frac{e^{-\abs{\vxp-\vx}^2/(4t)}}{4\pi t}\,d\vxp dt.\\
                                   &\leq\int_{0}^{\veps}\frac{e^{-d_B^2/(4t)}}{4\pi t}\int_{\Rd{2}}\sigma(\vxp)\delta(g(\vxp))(1-\Psi(\vxp))\,d\vxp dt,
          \end{split}
        \end{equation}
        since $(1-\Psi)$ is nonzero outside of $B$, i.e. on $\overline{\bdry\setminus B}$. The integral in $\vxp$ is bounded by the maximum value of $\sigma$ restricted to $\Gamma$ and the length $\Gamma$. Using the change of variables $a = d_B^2/(4t)$ lets us bound the integral in $t$,
        \begin{equation}
          \int_{0}^{\veps}\frac{e^{-d_B^2/(4t)}}{4\pi t}\, dt = \int_{d_B^2/(4\veps)}^{\infty}\frac{e^{-a}}{4\pi a}\, da\leq \frac{\veps}{\pi d_B^2}e^{-d_B^2/(4\veps)}.
        \end{equation}
        This lets us use the local graph parametrization  \eqref{eq:loc_param} and subsequent Taylor expansions in the following calculation without referring to the boundary outside this coordinate interval $U\cap \Gamma$.

Inside the localized integral we may replace the curve by the local graph \eqref{eq:loc_param},
\begin{equation}
  \label{eq:distcurveint3}
  \begin{split}
  \int_{\R^2}\Psi(\vxp)\sigma(\vxp)\delta(g(\vxp))\,d \vxp
  &=
    \int_\bdry\Psi(\vga(x_1'))\frac{\sigma(\vga(x_1'))}{\abs{\nabla g(\vga(x_1'))}}\abs{\vga'(x_1')}\,dx_1'\\
  &=  \int_\R\psi(x_1')\sigma(\vga(x_1'))\sqrt{1+(\curve'(x_1'))^2}\,dx_1',
      \end{split}
    \end{equation}
    since $g$ is the signed distance function, so $|\nabla g|=1$ on $\bdry$, and using properties of the Dirac distribution $\delta$. Here $\psi(x_1')=\Psi(\vga(x_1'))$ is the restriction of $\Psi$ to the local graph.

Recall that, in the local coordinate system, the target point $\vx = (0,r)$, where $r>0$ is the distance to the closest boundary point $\vb = (0,0)$. With the change of variables $z = \sqrt{t}$ and $\veta = z \vxi$, where $\veta = (\eta_1,\eta_2)$, we obtain
\begin{equation}\label{eq:slpl_rotated}
  \slpl[\sigma](\vx)  = \lim\limits_{\hat{\delta}\to 0^+} \frac{1}{(2\pi)^2}  \int_{\hat{\delta}}^{\sqrt{\veps}} \int_{\R^2} \int_{\R^2}  \Psi(\vxp)\sigma(\vxp)\delta(g(\vxp))  \, e^{-i \veta \cdot \vxp / z} \, d\vxp \;  e^{-|\veta|^2} \, e^{i \veta \cdot \vx / z}  \, \frac{1}{z^2} \, d\veta \; 2z \, dz.
\end{equation}
The limit of $\hat{\delta}$ to $0^{+}$ is a technicality needed for the following intermediate steps in order to avoid the case $z=0$. Using \eqref{eq:loc_param} with the scaled variable $x_1' = uz$, and applying \eqref{eq:distcurveint3} to \eqref{eq:slpl_rotated}, we obtain
\begin{equation}\label{eq:slp_rottransscaled}
  \begin{split}
    \slpl[\sigma](\vx) = \lim\limits_{\hat{\delta}\to 0^+}\frac{1}{2\pi^2} \int_{\hat{\delta}}^{\sqrt{\veps}} \int_{\R^2}e^{-|\veta|^2}\int_{\R}\psi(uz) \sigma(\vga(uz))\sqrt{1+(\gamma'(uz))^2} e^{-i\eta_1 u}e^{-i\eta_2 \curve(uz)/z}\,du \, e^{i\eta_2 r /z}\,d\veta \, dz
      \end{split}
\end{equation}
with an error of $\ordo{e}{-d_B^2/(4\veps)}$ due to neglecting contributions from $\bdry\setminus B$.

\subsubsection{Reduction in the $u$ and $\eta_1$ variables}\label{subsubsec:reduce_ueta}

Let us first expand the factors in the integrand of \cref{eq:slp_rottransscaled} in powers of $z$ using the local parameterization $x_1 = uz$. All quantities are expanded about the boundary point $\vb$, corresponding to $u=0$, such that for a smooth function $\varphi$ we have
\begin{equation}
  \varphi(uz) = \sum_{n=0}^{N}\frac{\varphi^{(n)}(0)}{n!}(uz)^n + R_{N+1}(uz),
\end{equation}
where $R_{N+1}(uz)$ is a smooth remainder. We now show that
\begin{equation}
 I_R = \lim\limits_{\hat{\delta}\to 0^+}\int_{\hat{\delta}}^{\sqrt{\veps}} \int_{\R}\int_{\R}\psi(uz) R_{N+1}(uz) e^{-\eta_1^2} e^{-i\eta_1 u}\,du \, d\eta_1 \, dz = \ordo{\veps}{(N+2)/2}.
\end{equation}
Assume that $\abs{\varphi^{(N+1)}/(N+1)!}$ is bounded on the support $U$ of $\psi$ by the constant $C_{\varphi}$, then by Taylor's theorem we have $\abs{R_{N+1}(uz)}<C_{\varphi}\abs{uz}^{N+1}$. Substituting this into $I_R$ and changing the order of integration to carry out the $\eta_1$ integral first gives
\begin{equation}\label{eq:I_R_bound}
    \abs{I_R} \leq C_{\varphi}\int_{0}^{\sqrt{\veps}}z^{N+1}\int_{\R}\abs{u}^{N+1}e^{-u^2/4}\,du \, dz \leq C_{R}\veps^{(N+2)/2},
  \end{equation}
where $C_R$ contains $C_{\varphi}$, the Gaussian moments and the other prefactors.

The layer density $\sigma$ satisfies
\begin{equation}
\sigma(\vga(uz)) = \sum_{n=0}^{N}\frac{\sigma_{b}^{(n)}}{n!}(uz)^n +R_{N+1}(uz),
\end{equation}
assuming $z<\sqrt{\veps}$ is small. Here $\sigma_b^{(n)} = \sigma^{(n)}(\vb)$ denotes the $n$th derivative of $\sigma$ with respect to the tangential parameterization at $\vb$.

Note that the boundary $\bdry$ given locally by  $\vga(x_1) = (x_1, \curve(x_1))$ satisfies
both $\curve(0) = 0$ and $\curve'(0) = 0$ since the tangent direction has been
aligned with the $x_1$-axis. Expanding $\curve(uz)$ gives
\begin{equation}
\curve(uz) = \frac{\kappa_{b}}{2}u^2z^2 + \sum_{n=3}^{N}\frac{\curve^{(n)}(0)}{n!}(uz)^n +R_{N+1}(uz),
\end{equation}
where $\kappa_{b} = \curve''(0)$ is the curvature at $\vb$. The oscillatory factor 
involving the boundary geometry can for a fixed $\eta_2$ be expanded as
\begin{equation}
e^{-i\eta_2\curve(uz)/z} = \sum_{n=0}^{N}\frac{(-i\eta_2)^n}{n!}\left(\frac{\curve(uz)}{z}\right)^n +u\, R_{N+1}(uz),
\end{equation}
 since $\curve(uz)/z = \ordo{z}{}$. We note that the error in a truncated Taylor expansion of $e^{-i\eta_2\curve(uz)/z}$ is a polynomial in $\eta_2$, assuming $uz$ constant. The integrand in $\slpl$ contains $e^{-\eta_2^2}$, thus the integral in $\eta_2$ is well defined and the error bound \eqref{eq:I_R_bound} holds. Finally, the Jacobian factor from the arc-length parameterization can be expanded as
\begin{equation}
\sqrt{1 + (\curve'(uz))^2} = 1+\sum_{n=1}^{N}(-1)^{n-1}\frac{(2n)!}{(2n-1)4^n(n!)^2}(\curve'(uz))^{2n} + R_{2(N+1)}(uz),
\end{equation}
since $\curve'(uz) = \ordo{z}{}$. As shown above in \cref{eq:I_R_bound}, truncation after total algebraic $z$-degree $N$ omits only terms of $\ordo{\veps}{(N+2)/2}$ from $\slpl$.

Collecting powers of $u$ and grouping Taylor series coefficients into functions $\alpha_n(\eta_2,z)$ of $\eta_2$ and $z$ gives 
\begin{align}
  \slpl[\sigma](\vx) &= \lim\limits_{\hat{\delta}\to 0^+}\frac{1}{2\pi^2}\int_{\hat{\delta}}^{\sqrt{\veps}} \int_{\R}e^{-\eta_2^2}e^{i\eta_2 r/z}\sum_{n=0}^{N}\alpha_n(\eta_2,z)\int_{\R}\int_{\R}\psi(uz) \, u^{n} e^{-\eta_1^2}e^{-i\eta_1 u}\,du \,\,d\eta_1 \,d\eta_2 \, dz + \ordo{\veps}{(N+2)/2}.
\end{align}
To reduce this expression to an integral involving only the variable $z$ we note that the integrand in $u$ is smooth and compactly supported for $z>0$, so we may switch the order of integration of $u$ with $\eta_1$. Let
\begin{equation}
      I^{\psi}_{n} := \int_\R\int_\R \psi(uz)u^ne^{-\eta_1^2}e^{-iu\eta_1} \,du\,d\eta_1,   
\end{equation}
and use that the Fourier transform of a Gaussian is again a Gaussian,
\begin{equation}\label{eq:FouriertranspairsGauss}
  \mathcal{F}\left[e^{-\eta_1^2}\right](u) = \sqrt{\pi}e^{-u^2/4},
\end{equation}
which gives
\begin{equation}
   I^{\psi}_{n} = \sqrt{\pi}\int_\R\psi(uz)u^ne^{-u^2/4}\,du.   
\end{equation}
Assume $\psi(uz) = 1$ for $\abs{u}\leq\varrho/z$ for some $\varrho>0$, then $1-\psi(uz) = 0$ unless $\abs{u}>\varrho/z$. The identity $\psi(uz) = 1 - (1-\psi(uz))$ allows us to write
\begin{equation}
  \begin{split}\label{eq:Ipsi}
    I^{\psi}_{n} &= \sqrt{\pi}\int_\R u^ne^{-u^2/4}(1-(1-\psi(uz)))\,du \\
                 &= \sqrt{\pi}\int_\R u^ne^{-u^2/4}\,du - \sqrt{\pi}\int_{\abs{u}>\frac{\varrho}{z}}(1-\psi(uz))u^ne^{-u^2/4}\,du.
  \end{split}
\end{equation}
Straightforward calculations show that $I^{\psi}_{n} = w_n + \ordo{e}{-1/z^2}$, where
\begin{equation}\label{eq:gauss_moment_identity}
   w_{n} := \sqrt{\pi}\int_\R u^ne^{-u^2/4}\,du =   \begin{cases}
   2\pi\dfrac{n!}{\left(\frac{n}{2}\right)!},
    & n \text{ even},\\
    0,
    & \text{otherwise.}
  \end{cases}
 \end{equation}

We now show that the second term in \cref{eq:Ipsi}, denoted $I_w$, is exponentially small in $z$. Let $C = \sqrt{\pi}\|1-\psi\|_{L^{\infty}}$ and split the Gaussian $e^{-u^2/4}=e^{-u^2/8}e^{-u^2/8}$, then
\begin{equation}
  \begin{split}\label{eq:Iw_est}
    \abs{I_{w}} \leq C  \int_{\abs{u}>\frac{\varrho}{z}}|u|^ne^{-u^2/4}\,du \leq C e^{-\varrho^2/(z^2 8)} \int_{\R}|u|^ne^{-u^2/8}\,du = \ordo{e}{-C_{\varrho}^2/z^2},
  \end{split}
\end{equation}
with constant $C_{\varrho} < \varrho/8$. Thus the integral of $z^m I_w$ in $z$, for each fixed integer $m$, is $\ordo{e}{-C_{\varrho}^2/\veps}$. We can now safely pass the limit $\hat{\delta}\to 0^+$, and the expression for $\slpl$ can now be reduced to integrals in $\eta_2$ and $z$,
\begin{equation}\label{eq:slplinetaz}
  \slpl[\sigma](\vx) = \frac{1}{\pi}\int_{0}^{\sqrt{\veps}} \sum_{n=0}^{\lfloor N/2 \rfloor}\frac{(2n)!}{n!}\int_{\R}e^{-\eta_2^2}\alpha_{2n}(\eta_2,z) \, e^{i\eta_2 r/z}\,d\eta_2 \, dz + \ordo{\veps}{(N+2)/2},
\end{equation}
with leading order errors from the Taylor expansions. Next, we consider integration with respect to the variable $\eta_2$.
\subsubsection{Reduction in the $\eta_2$ variable}
We expand the $\eta_2$-dependent factors in \cref{eq:slplinetaz} in powers of $\eta_2$, and we obtain $\sum_{n = 0}^N \beta_n(z)\eta_2^n$, where $\beta_n(z)$ consists of the contributions from the Taylor expansions of the density and the boundary geometry, combined with the constants $w_n$ from \cref{eq:gauss_moment_identity}. For $n\geq 0$, we define the
$n$th \emph{Fourier--Gaussian moment} by
\begin{equation}\label{eq:I_n}
    I_{n}(c):= \int_\R \eta^n e^{-\eta^2} e^{i\eta c}\,d\eta,\quad c\in\R.
\end{equation}
For each $n$, the function $I_n$ admits the explicit representation
\begin{equation}\label{eq:I_n_exp}
  I_{n}(c) = \sqrt{\pi}\frac{i^n}{2^{n}}H_n\left(\frac{c}{2}\right)e^{-\frac{c^2}{4}},
\end{equation}
where $H_n(c)$ denotes the physicists' Hermite polynomial of degree $n$. Substituting this into \cref{eq:slplinetaz}, we obtain
\begin{equation}
    \slpl[\sigma](\vx) = \frac{1}{\pi}\int_{0}^{\sqrt{\veps}} \sum_{n=0}^{N}\beta_{n}(z)
    I_{n}\!\left(\frac{r}{z}\right) \, dz = \frac{1}{\sqrt{\pi}}\int_{0}^{\sqrt{\veps}} \sum_{n=0}^{N}\beta_{n}(z)
    \frac{i^n}{2^{n}}H_n\!\left(\frac{r}{2z}\right)e^{-\frac{r^2}{4z^2}} \, dz + \ordo{\veps}{(N+2)/2},
\end{equation}
where $r$ is the distance between the target point and its closest point on the boundary.

\subsubsection{Reduction in the $z$ variable}\label{subsubsec:slppois_z}
Expanding each $\beta_n(z)$ in powers of $z$ and collecting terms
with respect to powers of $z$, we obtain
\begin{equation}
  \begin{split}
  \slpl[\sigma](\vx)
  = \sum_{m=0}^{N}\sum_{n=0}^{N}\zeta_{m,n}\frac{i^n}{2^{n}}
  \int_{0}^{\sqrt{\veps}}H_n\!\left(\frac{r}{2z}\right)
  e^{-\frac{r^2}{4z^2}}z^{m}\,dz + \ordo{\veps}{(N+2)/2},
  \end{split}
\end{equation}
where the coefficients $\zeta_{m,n}$ contain all prefactors together with the
contributions from the density $\sigma$ and the local geometry of $\bdry$.
We now make the change of variables $\omega = r/(2z)$, which maps
$(0,\sqrt{\veps})$ to $(r/(2\sqrt{\veps}),\infty)$, and use the expansion
$H_n(\omega) = \sum_{k=0}^n h_{n,k}\omega^k$ where $h_{n,k}$ are the coefficients of the Hermite polynomials, and obtain
\begin{equation}
 \begin{split}
    \slpl[\sigma](\vx) 
                        &= \sum_{m=0}^{N}\sum_{n=0}^{N}\zeta_{m,n}\frac{i^n}{2^{n}}\sum_{k=0}^{n}h_{n,k}\int_{r/(2\sqrt{\veps})}^{\infty}\omega^k
                          e^{-\omega^2}\frac{r^{m+1}}{2^{m + 1}}\frac{1}{\omega^{m+2}}\,d\omega\\
                       &= \sum_{m=0}^{N}\sum_{n=0}^{N}\zeta_{m,n}\frac{i^nr^{m+1}}{2^{m + n + 1}}\sum_{k=0}^{n}h_{n,k}Q_{m+2-k}\left(\frac{r}{2\sqrt{\veps}}\right) +\ordo{\veps}{(N+2)/2},
      \end{split}
    \end{equation}
where we define
\begin{equation}\label{eq:int_inc_gamma_def}
  Q_p(c) := \int_c^\infty\frac{e^{-x^2}}{x^p}\,dx, \quad c>0,\,p\in\mathbb{Z}.
\end{equation}
Equivalently,  $Q_p(c) = \frac{1}{2}\Gamma\left(\frac{1-p}{2},c^2\right)$ where $\Gamma(\cdot,\cdot)$ is the upper incomplete Gamma function. 
    
Introducing the scaled variable $c = r/\sqrt{\veps}$, we obtain
\begin{equation}\label{eq:slp_closedsum}
  \begin{split}
    \slpl[\sigma](\vx) = \sum_{m=0}^{N} \veps^{(m+1)/2} \left(\frac{c}{2}\right)^{m+1} \sum_{n=0}^{N}\zeta_{m,n}\frac{i^n}{2^{n}}\sum_{k=0}^{n}h_{n,k}Q_{m+2-k}\left(\frac{c}{2}\right) + \ordo{\veps}{(N+2)/2},
      \end{split}
    \end{equation}
    which establishes the desired expression.

 The functions $Q_p$ are exponentially small for large argument. For example, $Q_0(6)=\frac{\sqrt{\pi}}{2}\erfc(6) \approx 1.34\cdot 10^{-17}$, so in double precision, terms involving $Q_p(c/2)$ are negligible once $c \gtrsim 12$, up to polynomial prefactors. This motivates why it is natural to keep $c=r/\sqrt{\veps}$ fixed when writing the local expansion in powers of $\sqrt{\veps}$. The coefficients are functions of the scaled distance $c$, while the explicit powers of $\sqrt{\veps}$ determine the asymptotic order. Thus the formula is an expansion in $\veps$, with $c$ treated as an $\ordo{1}{}$ parameter.
    
The representation in terms of $Q_p(c/2)$ is written for $r>0$, equivalently $c=r/\sqrt{\epsilon}>0$, since $Q_p$ was defined with a positive lower integration limit. This does not exclude boundary targets. If $r=0$, one should instead use the preceding representation in terms of the Fourier--Gaussian moments $I_n(r/z)$, which are well defined at $r=0$. However, after the $Q_p$ terms have been combined into the final coefficient formulas, the same boundary values are obtained by taking the limit $c\to 0^+$.

\subsubsection{Obtaining $\slpl$ to a given order}\label{subsubsec:slppois_order}
To compute $\slpl$ from \eqref{eq:slp_closedsum}, we require $Q_p$, $p=-2,\ldots,N+2$, which are readily obtained from the recursion formula for the incomplete Gamma function $\Gamma$ in \cite{NISTDLMF}:
\begin{equation}\label{eq:int_inc_gamma_recur}
  Q_{p-2}(c) = \frac{1-p}{2}Q_{p}(c) + \frac{c^{1-p}}{2}e^{-c^2} \quad \text{for } c > 0,
\end{equation}
for any integer $p$. Since the recursion \cref{eq:int_inc_gamma_recur} steps by two, even and odd values of $p$ form separate chains. We use $Q_0(c)$ as the base case for even $p$, and $Q_1(c)$ as the base case for odd $p$. For $p = -2,-1,0,1$ we have that
\begin{equation}
  Q_{-2}(c) = \frac{c}{2}e^{-c^2}+\frac{\sqrt{\pi}}{4}\erfc(c),\quad  Q_{-1}(c) = \frac{1}{2}e^{-c^2},\quad  Q_{0}(c) = \frac{\sqrt{\pi}}{2}\erfc(c),\quad Q_{1}(c) = \frac{1}{2}E_1(c^2).
\end{equation}
The coefficients $\zeta_{m,n}$ require some tedious bookkeeping, starting with 
the coefficients $\alpha_n$ based on the Taylor series expansions of the layer density and geometric quantities.

To derive the expansion of $\slpl$ to $\ordo{\veps}{(P+1)/2}$ for a given $P$, we must set $N = P-1$. We verify \cref{eq:slp_closedsum} for $P=4$ against the closed-form expression derived in \cite{lghtwght_ps2026}. The coefficients in \eqref{eq:slplgeneral} for a given target point $\vx$ are
\begin{equation}
  A_p = \left(\frac{c}{2}\right)^{p} \sum_{n=0}^{N}\zeta_{p-1,n}\frac{i^n}{2^{n}}\sum_{k=0}^{n}h_{n,k}Q_{p+1-k}\left(\frac{c}{2}\right),
\end{equation}
where $A_0 = 0$ and
\begin{align}\label{eq:slplocO5full}
  A_1 &= \frac{\sigma_{b}}{2}\left(\frac{2 e^{-\frac{c^2}{4}}}{\sqrt{\pi }} - c \,\mathrm{erfc}\left(\frac{c}{2}\right)\right)\\
  A_2 &= \frac{c  \kappa_{b} \sigma_{b} }{4}\left(\frac{2 e^{-\frac{c^2}{4}}}{\sqrt{\pi }}-c \,\mathrm{erfc}\left(\frac{c}{2}\right)\right)\\
    A_3 &= \frac{\sigma_{b} \kappa_{b}^2}{12} \left(\frac{\left(4 c^2+1\right) e^{-\frac{c^2}{4}}}{\sqrt{\pi }}-2 c^3 \mathrm{erfc}\left(\frac{c}{2}\right)\right) +  \frac{\sigma^{\prime\prime}_{b}}{12}\left(\frac{2 \left(2-c^2\right) e^{-\frac{c^2}{4}}}{\sqrt{\pi }} + c^3 \mathrm{erfc}\left(\frac{c}{2}\right)\right)\\
  A_4 &=  \left[
\frac{c  \kappa_{b}^3 \sigma_{b} }{16} \left(\frac{\left(6 c^2-2\right) e^{-\frac{c^2}{4}}}{\sqrt{\pi }} - 3 c^3 \mathrm{erfc}\left(\frac{c}{2}\right)\right) + \frac{c  \kappa_{b} \sigma^{\prime\prime}_{b}}{8} \left(\frac{2 \left(2-c^2\right) e^{-\frac{c^2}{4}}}{\sqrt{\pi }}+c^3 \mathrm{erfc}\left(\frac{c}{2}\right)\right) \right]\notag \\
   & + 
\left[ 
\frac{c \sigma_{b} \gamma_{b}^{(4)}}{48}\!\left(\frac{2 \left(2-c^2\right) e^{-\frac{c^2}{4}}}{\sqrt{\pi }} + c^3 \mathrm{erfc}\left(\frac{c}{2}\right)\right) + \frac{c \sigma^{\prime}_{b} \gamma_{b}^{(3)}}{12}\!\left(\frac{(4-2 c^2)e^{-\frac{c^2}{4}}}{\sqrt{\pi}} + c^3 \mathrm{erfc}\left(\frac{c}{2}\right)\right) \right],
\end{align}
where $\curve^{(n)}$ is the $n$th derivative of the graph $\curve$ in the local tangent coordinates. This completes the derivation of the asymptotic expansion \cref{eq:slp_closedsum} for the local part $\slpl[\sigma](\vx)$.
\subsection{The double layer potential $\dlp$}\label{subsec:dlppois}
 We define the local double layer contribution by $\dlpl[\mu]:=\ul[F]$, where the dipole source is
\begin{equation}
F(\vx)=-\nabla\cdot\bigl(\vnu(\vx)\mu(\vx)\delta(g(\vx))\bigr),
\end{equation}
with $\vnu(x_1)$ the outward directed unit normal. Passing to Fourier space, the spatial derivative becomes multiplication by $i\vk$, so that
\begin{equation}
\calF\!\left[-\nabla\cdot (\vnu\mu\delta(g))\right]
=
-i\vxi\cdot \calF[\vnu\,\mu\,\delta].
\end{equation}
After the same localization and change of variables as in \cref{subsubsec:slppois_locparam}, the normal satisfies $\vnu(x_1)=\frac{(\curve'(x_1),-1)}{\sqrt{1+(\curve'(x_1))^2}}$, and the normalization factor $\sqrt{1+(\curve'(x_1))^2}$ in the denominator cancels the Jacobian from the arc-length parameterization. Let $\tilde{\vnu}(uz) = (-\curve'(uz),1)$, then
\begin{equation}\label{eq:dlpl_rotated}
  \dlpl[\mu](\vx)
  =
  \lim\limits_{\hat{\delta}\to 0^{+}}
  \frac{1}{2\pi^{2}}
  \int_{\hat{\delta}}^{\sqrt{\veps}}
  \int_{\R^2} e^{-|\veta|^2}
  \int_{\R}\psi(uz)
  i\,\veta\cdot\tilde{\vnu}(uz)\,\mu(uz)\,
  e^{-i\eta_1 u}e^{-i\eta_2 \curve(uz)/z}\,du\,
  e^{i\veta\cdot \vx/z}\,d\veta\,
  \frac{1}{z}\,dz,
\end{equation}
with an exponentially small error $\ordo{e}{-d_B^2/(4\veps)}$ due to contributions from the neglected part of the curve. The additional factor of $\veta$ modifies only the algebraic structure of the intermediate expansions. The subsequent reduction is identical to that of \cref{subsec:slppois}, and is therefore omitted. We simply state the resulting expansions of the form
\begin{equation}\label{eq:dlplgeneral}
  \dlpl[\mu] (\bm{x}) = \sum_{p=0}^P\veps^{p/2}A_p + \ordo{\veps}{(P+1)/2}.
\end{equation}
  For $P=4$, we have the coefficients
\begin{align}\label{eq:dlplocO5full_012}
A_0 &= -\frac{\mu_{b}}{2} \mathrm{erfc}\left(\frac{c}{2}\right),\quad A_1 = -\frac{\kappa_{b} \mu_{b}}{2 \sqrt{\pi }}e^{-\frac{c^2}{4}},\quad A_2= -\left[ 
\frac{3 c \kappa_{b}^2 \mu_{b} }{8 \sqrt{\pi }}e^{-\frac{c^2}{4}} 
+ \frac{c  \mu^{\prime\prime}_{b}}{4}  \left(\frac{2 e^{-\frac{c^2}{4}}}{\sqrt{\pi }}-c\, \mathrm{erfc}\left(\frac{c}{2}\right)\right) \right], \\
  A_3 &= -\frac{5 \left(c^2-2\right) \kappa_{b}^{3} \mu_{b}}{16 \sqrt{\pi }}e^{-\frac{c^2}{4}} - \frac{\gamma^{(4)}_{b}\mu_{b}}{4 \sqrt{\pi }}e^{-\frac{c^2}{4}} - \frac{ \kappa_{b} \mu^{\prime\prime}_{b}}{4}\left(\frac{2 \left(c^2+1\right) e^{-\frac{c^2}{4}}}{\sqrt{\pi }}-c^3 \mathrm{erfc}\left(\frac{c}{2}\right)\right) \notag\\\label{eq:dlplocO5full_3}
  &- \frac{\gamma^{(3)}_{b} \mu^{\prime}_{b}}{12} \left(\frac{2e^{-\frac{c^2}{4}} \left(c^2+4\right)}{ \sqrt{\pi }}-c^3 \mathrm{erfc}\left(\frac{c}{2}\right)\right) \\
A_4 &=  -\left[ 
\frac{35 c  \kappa_{b}^4 \mu_{b} \left(c^2-6\right) }{128 \sqrt{\pi }}e^{-\frac{c^2}{4}}  
+ \frac{5 c\left(\gamma^{(3)}_{b}\right)^{2} \mu_{b} }{12 \sqrt{\pi }} e^{-\frac{c^2}{4}} 
+ \frac{5 c \kappa_{b}\gamma^{(4)}_{b} \mu_{b}}{8 \sqrt{\pi }}e^{-\frac{c^2}{4}}
\right]   \notag\\
& - 5c \left[ 
\frac{\kappa_{b}^2 \mu^{\prime\prime}_{b}}{16} \left(\frac{2 \left(c^{2}+1\right) e^{-\frac{c^2}{4}}}{\sqrt{\pi }}-c^{3} \mathrm{erfc}\left(\frac{c}{2}\right)\right)
+
 \frac{\kappa_{b} \gamma^{(3)}_{b} \mu^{\prime}_{b} }{24 }\left(\frac{2 \left(c^2+4\right)e^{-\frac{c^2}{4}}}{\sqrt{\pi }} - c^3  \mathrm{erfc}\left(\frac{c}{2}\right) \right) \right] \notag\\ \label{eq:dlplocO5full_4}
&-\frac{c  \mu^{(4)}_{b}}{48}\left(c^3 \mathrm{erfc}\left(\frac{c}{2}\right)-\frac{2 \left(c^2-2\right) e^{-\frac{c^2}{4}}}{\sqrt{\pi }}\right),
\end{align}
which is in agreement with the corresponding expansion for $\dlpl$ in \cite{lghtwght_ps2026}.

The coefficients in \eqref{eq:dlplocO5full_012}--\eqref{eq:dlplocO5full_4} are written for the interior side, corresponding to $\sgn=-1$ in \eqref{eq:sgn_def} If the same coefficient profiles are used for targets on either side of $\bdry$, then the even-order terms in the double layer expansion must be multiplied by $-\sgn$, while the odd-order terms are unchanged. Thus
\begin{equation}
  \dlpl[\mu](\vx)  =  \sum_{\substack{0\leq p\leq P\\ p\ \mathrm{even}}}
  \veps^{p/2}(-\sgn)A_p
  +
  \sum_{\substack{0\leq p\leq P\\ p\ \mathrm{odd}}}
  \veps^{p/2}A_p
  + \ordo{\varepsilon}{(P+1)/2}.
\end{equation}
In particular, the leading term becomes
\begin{equation}
  A_0 =
  \frac{\rho\mu_b}{2}\operatorname{erfc}\!\left(\frac{c}{2}\right),
\end{equation}
which makes the double layer jump explicit.

\subsection{The Volume Potential $\vp$}\label{subsec:vppois}
Deriving asymptotic expansions for the local part of the volume potential $\vp$ follows the same framework as in \cref{subsec:slppois}, and we therefore use the same notation and scaling throughout. We define $\vpl[f](\vx) := \ul[F](\vx)$ with $F = f\ind$, where $\ind$ is the indicator function of $\dom$, and $f$ a smooth source density on $\dom$. The essential difference from the layer potentials is that the source is supported on a two-dimensional region rather than on a curve, which will introduce a Heaviside factor in the inward directed normal direction.

We place the target point at $\vx$ at $(0,r)$ and let the closest boundary point $\vb$ be located at the origin $(0,0)$. As for $\slpl$, we assume $r>0$ and then z-integrand in \eqref{eq:vppois} is continuous on $(0,\sqrt{\veps}]$.

For $\veps$ sufficiently small, we restrict attention to a neighborhood in which the boundary admits the same local parameterization as in \cref{subsec:slppois}, namely $\vga(\vxp) = (x'_1,\curve(x'_1))$, with exponentially small error. In this representation, points in $\dom$ satisfy $x'_2 > \curve(x'_1)$, and we write the indicator function as $\ind(\vxp) = \theta(x'_2 - \curve(x'_1))$, where $\theta$ is the Heaviside function. 

Let $\Psi$ be separable such that $\Psi(\vx) = \psi_1(x_1')\psi_1(x_2')$ and introduce the variables $t = z^2$, $\veta = z\vxi$, and $\vu = \vxp/z$. Then the local part takes the form
\begin{equation}\label{eq:vppois}
  \vpl[f](\vx)
  = \lim\limits_{\hat{\delta}\to 0^+}\frac{1}{2\pi^2}\int_{\hat{\delta}}^{\sqrt{\veps}} z \int_{\R^2} e^{-|\veta|^2}
  \int_{\R^2} \Psi(\vu z)f(\vu z)\, \theta\!\left(u_2 - \frac{\curve(u_1 z)}{z}\right)
  e^{-i\veta\cdot\vu}e^{i\veta\cdot\vx/z}\, d\vu\, d\veta\, dz,
\end{equation}
with an error of $\ordo{e}{-d_B^2/(4\veps)}$ as for $\slpl$. It follows from repeating the steps in \cref{eq:slp_err_local}:
\begin{equation}
\vpl[(1-\Psi)f](\vx) \leq\lim\limits_{\hat{\delta}\to 0^+} \int_{\hat{\delta}}^{\veps}\frac{e^{-d_B^2/(4t)}}{4\pi t}\int_{\Omega}f(\vxp)(1-\Psi(\vxp))\,d\vxp dt.
\end{equation}

All quantities in \cref{eq:vppois} are expanded in powers of $z$ about the target point $\vx = (0,r)$. For multivariate expansions we use total degree, so that
\begin{equation}\label{eq:fun_exp}
  f(\vu z)
  = \sum_{m+n\le N} \frac{f^{(m,n)}}{m!\,n!} u_1^m u_2^n z^{m+n}
  +R_{N+1}(uz),
\end{equation}
where $f^{(m,n)} = \partial_1^m\partial_2^n f(\vx)$. The geometry enters only through the expansion of $\curvex(u_1,z) = \curve(u_1 z)/z$, which we write as
\begin{equation}\label{eq:curve_tail}
\curvex(u_1,z)
  = \frac{\kappa_b}{2} u_1^2 z
  + \sum_{n=3}^{N+1} \frac{\curve^{(n)}(0)}{n!} u_1^n z^{n-1}
  +uR_{N+1}(uz).
\end{equation}
For all $R_{N+1}(uz)$ above the bound \eqref{eq:I_R_bound} holds for $\vpl$ as well, which follows from the exponential decay of the Gaussian.

With these expansions in place, the reduction proceeds in the same sequence as for the single layer potential. The Heaviside factor affects only the treatment of the $u_2$-integration, while the subsequent steps: reduction in $\eta_1$ and $u_1$, followed by $\eta_2$, and finally the $z$-integration, are identical to those in \cref{subsec:slppois}. This leads to an asymptotic expansion of $\vpl[f](\vx)$ in powers of $\veps^{p/2}$, with coefficients determined by derivatives of $f$ and the graph $\curve$ at $\vb$.

\subsubsection{Reduction in the $u_2$ variable}\label{subsubsec:vppois_du2}
The $u_2$-dependent part of the integrand in \cref{eq:vppois} is $\psi_2(u_2z)u_2^n \theta\!\left(u_2 - \curvex(u_1,z)\right)e^{-i\eta_2 u_2}$, with $n \le N$. We proceed as in \cref{subsec:slppois}, which requires a distributional treatment. First, the cutoff function $\psi_2$ makes the integral in $u_2$ well-defined and, as in the reduction leading to \cref{eq:Iw_est}, it may be removed with an exponentially small error in $\veps$.

Using the Fourier transform of the Heaviside function $\calF[\theta](\xi) = \pi\delta(\xi) + \pv\frac{1}{i\xi}$, where $\pv$ denotes the Cauchy principal value, together with the translation property
\begin{equation}\label{eq:Fouriertranspairs1}
  \calF[f(x-a)] = e^{-i a \xi}\hat f(\xi),
\end{equation}
and differentiation property
\begin{equation}\label{eq:Fouriertranspairs2}
  \mathcal{F}\left[u^{n}\right](\eta_1) = 2\pi i^{n}\delta^{(n)}(\eta_1),\quad n\geq 0,
\end{equation}
 we obtain
\begin{equation}\label{eq:fourier12heaviside}
\calF\left[u_2^{n} \theta\!\left(u_2 - \curvex(u_1,z)\right)\right](\eta_{2})
= i^n \frac{\partial^{n}}{\partial \eta_{2}^{n}}
\left(e^{-i\eta_{2}\curvex(u_1,z)}\left(\pi\delta(\eta_{2}) + \pv\frac{1}{i\eta_2}\right)\right).
\end{equation}
Substituting \cref{eq:fun_exp} and \cref{eq:fourier12heaviside} into \cref{eq:vppois} and carrying out the $u_2$-integration, we obtain
\begin{equation}\label{eq:vppois_int1}
\begin{split}
\vpl[f](\vx)
&= \lim\limits_{\hat{\delta}\to 0^+}\frac{1}{(2\pi)^2}\int_{\hat{\delta}}^{\sqrt{\veps}}\int_{\R^2} e^{-|\veta|^2}
\sum_{m+n\le N} \frac{f^{(m,n)}_0}{m!n!} z^{m+n+1}
\int_{\R} \psi_1(u_1z)u_1^m e^{-i\eta_1 u_1} e^{i\eta_2 r/z} e^{-i\eta_2 \curvex(u_1,z)} \\
&\quad  \sum_{k=0}^n \binom{n}{k}
\left(\curvex(u_1,z)\right)^{n-k}
i^k \partial_{\eta_2}^k \left(\pi\delta(\eta_2) + \pv\frac{1}{i\eta_2}\right)
\,du_1\, d\veta\, dz\, +\, \ordo{\veps}{(N+2)/2},
\end{split}
\end{equation}
where $\curvex(u_1,z)$ is expanded as in \cref{eq:curve_tail}.

The expression in \cref{eq:vppois_int1} separates naturally into a contribution from the $\delta$-term and a contribution from the principal value term. These have different algebraic structures and will be treated independently, but first we carry out the reduction steps in $(u_1,\eta_1)$ as in \cref{subsec:slppois}.

\subsubsection{Reduction in the $u_1$ and $\eta_1$ variables}\label{subsubsec:vppois_deta1_du1}
Every term in \cref{eq:vppois_int1} contains a factor $\psi_1(uz)u_1^m e^{-\eta_1^2}e^{-i\eta_1 u_1}$ for $m\geq 0$, that comes from the expansion of the source density $f$, the oscillatory factor $e^{-i\eta_2\curvex(u_1 z)}$, and powers of $\curvex(u_1,z)$. As in the single layer case, the integration in $(\eta_1,u_1)$ reduces to the moments  $w_m$ from \cref{eq:gauss_moment_identity} with an exponentially small error in $\veps$. All terms with odd $m$ vanish, while for even $m$ we obtain the prefactors $w_m = 2\pi\tfrac{m!}{(m/2)!}$.

\subsubsection{Reduction in the $\eta_2$ variable}\label{subsubsec:vppois_deta2}
After reduction in the variables $u_2$, $\eta_1$ and $u_1$, the expression \eqref{eq:vppois_int1} is
\begin{equation}\label{eq:vppois_int_eta2}
\begin{split}
\vpl[f](\vx)
&= \frac{1}{(2\pi)^2}\int_{0}^{\sqrt{\veps}} \sum_{m+n\leq N} z^{m+n+1}
\sum_{k=0}^n \binom{n}{k} i^k \\
&\quad  \int_{\R} e^{-\eta_2^2} e^{i\eta_2 r/z}
\alpha_{m,n,k}(z,\eta_2)
\partial_{\eta_2}^k\!\left(\pi\delta(\eta_2) + \pv\frac{1}{i\eta_2}\right)
\,d\eta_2\, dz \, + \ordo{\veps}{(N+2)/2},
\end{split}
\end{equation}
where $\alpha_{m,n,k}(z,\eta_2)$ collects the algebraic contributions from the Taylor expansion of $f$, the expansion of $\curvex(u_1,z)$ in \eqref{eq:curve_tail}, the expansion of $e^{-i\eta_2\curvex(u_1,z)}$, and the moments $w_{m}$.

For fixed $m,n,k$, the $\eta_2$-integrals in \cref{eq:vppois_int_eta2} take the form
\begin{equation}
\int_{\R}\phi_p(\eta_2)\,\delta^{(k)}(\eta_{2})\,d\eta_2,
\qquad
\int_{\R}\phi_p(\eta_2)\,\partial^k_{\eta_2}\left(\pv\frac{1}{i\eta_2}\right)\,d\eta_2,
\end{equation}
where
\begin{equation}\label{eq:phip}
  \phi_p(\eta_2) = e^{-\eta_2^2} e^{i\eta_2 r/z}\eta_2^{p},\quad  p \ge 0.
\end{equation}
Differentiation of the distributions produces, in each case, polynomials in the variable $\eta_2$. As a consequence, the integrals reduce to combinations of Fourier--Gaussian moments $I_p(r/z)$. Carrying out these reductions yields expressions that can be written as polynomials in $z$, with coefficients determined by the functions $\alpha_{m,n,k}(z,\eta_2)$, which collect the contributions from the Taylor expansions and the geometric terms.

The local part $\vpl$ therefore decomposes naturally into two distinct contributions: one arising from the $\delta$-term and one from the principal value term. We denote these by $\vpldd$ and $\vplpv$, respectively. Since these contributions have different algebraic structures, they are treated separately. We begin with the contribution associated with the $\delta$-term.
\subsubsection{Contribution from the $\delta$--term}

Using the distributional identity
\begin{equation}
\int_\R \phi_p(\eta_2)\,\delta^{(k)}(\eta_2)\,d\eta_2 = (-1)^k \phi_p^{(k)}(0),
\end{equation}
we evaluate the $\delta$-contribution by differentiating $\phi_p$ in \cref{eq:phip}. 
Writing $\phi_p$ as a Taylor expansion and differentiating termwise gives
\begin{equation}\label{eq:int_eta2_diracdiff}
\phi_p^{(k)}(\eta_2)\vert_{\eta_2=0}
= \frac{\partial^k}{\partial \eta_2^k}\left.\left(\sum_{j,\ell \ge 0}
\frac{(ir/z)^j}{j!}\frac{(-1)^\ell}{\ell!}
\eta_2^{p + j + 2\ell}\right)\right|_{\eta_2=0}
= \sum_{\substack{j,\ell \ge 0 \\ p + j + 2\ell = k}}
\frac{(ir/z)^j}{j!}\frac{(-1)^\ell}{\ell!}k!
\end{equation}
Substituting this into \cref{eq:vppois_int_eta2} yields
\begin{equation}\label{eq:vppois_int5_delta}
\vpldd[f](\vx)
= \lim\limits_{\hat{\delta}\to 0^+}\sum_{m+n\leq N} \sum_{k=0}^n \binom{n}{k} i^k
\int_{\hat{\delta}}^{\sqrt{\veps}} z^{m+n+1}\,\beta_{m,n,k}(z)\,dz,
\end{equation}
where $\beta_{m,n,k}(z)$ collects the resulting coefficients, consisting of geometric data, in particular $r$ and $\kappa_b$, Taylor coefficients of $f$, index-dependent constants arising from evaluation at $\eta_2 = 0$, the moments $w_{m}$, and powers of $z$ and $r/z$.

We investigate if there are any negative powers of $z$ in the integrand of $\vpldd$. By \cref{eq:vppois_int5_delta}, all terms carry a factor $z^{m+n+1}$, and for fixed $k$ the factor $(\curvex(u_1,z))^{\,n-k}$ contributes at most $z^{k-n}$. From \cref{eq:int_eta2_diracdiff}, the derivatives contribute at most $(r/z)^j$ with $0 \le j \le k$. Let $P(\eta_2)$ be a polynomial in $\eta_2$, then
\begin{align}\label{eq:zpower_count_delta}
&z^{m+n+1} (\curvex(u_1,z))^{\,n-k}\,
\partial_{\eta_2}^{k}\!\left(P(\eta_2)e^{-\eta_2^2}e^{-i\eta_2\curvex(u_1,z)}\right)\big|_{\eta_2=0} \\
&= \ordo{z}{m+n+1}\,\ordo{z}{k-n}\,\ordo{z}{-j}
= \ordo{z}{m+1+k-j}, \qquad 0 \le j \le k.
\end{align}
The lowest-order contribution occurs for $j = k$, yielding $\ordo{z}{m+1}$. Since $m \ge 0$, no terms of order $z^{-1}$ arise, so we can safely pass the limit  $\hat{\delta}$ to $0^+$, giving
\begin{equation}\label{eq:vppois_int6_delta}
\vpldd[f](\vx)
= \sum_{m=0}^{N} \zeta^\delta_{m} \int_0^{\sqrt{\veps}} z^{m+1}\,dz,
\end{equation}
where the coefficients $\zeta^\delta_m$ collect contributions from the Taylor coefficients of $f$, binomial factors, the moments $w_{m}$, the local geometric coefficients from the expansion of $\curvex(u_1,z)$, and the factors of $r$ arising from the evaluation of \cref{eq:int_eta2_diracdiff}.

\subsubsection{Contribution from the $\pv\frac{1}{i\eta_2}$--term}

Derivatives of $\pv \tfrac{1}{i\eta_2}$ are understood in the distributional sense, so for any $\phi \in \schwartz(\R)$,
\begin{equation}
  \int_\R \phi(\eta_2)\,\partial_{\eta_2}^{k}\!\left(\pv \frac{1}{i\eta_2}\right)\,d\eta_2
  = (-1)^k\, \pv\int_\R \phi^{(k)}(\eta_2)\, \frac{1}{i\eta_2}\,d\eta_2,
\end{equation}
since all derivatives of $\phi$ vanish at infinity. It follows that the derivatives may be transferred onto the smooth factor, and the resulting principal value integrals are evaluated using
\begin{equation}
  \pv\int_\R \frac{\phi(\eta_2)}{i\eta_2}\,d\eta_2  = \frac{1}{i}\int_0^{\infty}\frac{\phi(\eta_2)-\phi(-\eta_2)}{\eta_2}\,d\eta_2.
\end{equation}
We have $\phi(\eta_2) = e^{-\eta_2^2} e^{i\eta_2 r/z}\eta_2^{p}$ from \eqref{eq:phip} with $p\geq 0$ and write $e^{i\eta_2 r/z} = \cos(\eta_2 r/z) + i\sin(\eta_2 r/z)$. The difference $\phi(\eta_2)-\phi(-\eta_2)$ keeps only the odd part of $\phi$, and the odd part contains a factor of $\eta_2^p$ for $p>0$, so the quotient by $\eta_2$ is regular at the origin. Hence, for $p>0$, the resulting integrals are Fourier--Gaussian moments $I_{p-1}(r/z)$ from \eqref{eq:I_n}. For $p = 0$ we have $\sin(\eta_2 r/z)/\eta_2$ which is regular at $\eta_2 = 0$, and we have
\begin{equation}
  \pv\int_\R \frac{e^{-\eta_2^2} e^{i\eta_2 r/z}}{i\eta_2}\,d\eta_2 = \pi\erf\left(\frac{r}{2z}\right).
\end{equation}
Consequently,
\begin{equation}\label{eq:vppois_int8_pv}
  \vplpv[f](\vx)
  =
  \sum_{\substack{m < N \\ n < N}}
  \hat{\zeta}^{\pv}_{m,n}
  \lim\limits_{\hat{\delta}\to 0^+}\int_{\hat{\delta}}^{\sqrt{\veps}} z^{n+1} I_m(r/z)\,dz +
  \sum_{n < N}
  \tilde{\zeta}^{\pv}_{n}
  \lim\limits_{\hat{\delta}\to 0^+}\int_{\hat{\delta}}^{\sqrt{\veps}} z^{n+1} \erf\left(\frac{r}{2z}\right)\,dz,
\end{equation}
where $\hat{\zeta}^{\pv}_{m,n}$ and $\tilde{\zeta}^{\pv}_{n}$ collect the moments $w_m$, Taylor coefficients of $f$, and geometric contributions from $\curvex$.

The limit $\hat{\delta}$ to $0^+$ can now be passed. For fixed $r\geq 0$, the functions $I_m(r/z)$ and $\operatorname{erf}(r/(2z))$ remain bounded as $z\to 0^+$. For $r>0$, $I_m(r/z)$ is a polynomial in $r/z$ times $e^{-r^2/(4z^2)}$, and for $r=0$ it is simply $I_m(0)$. Hence each integrand in \eqref{eq:vppois_int8_pv} is bounded by $C z^{n+1}$ near $z=0$, which is integrable, so the limit $\hat{\delta}\to0^+$ is an ordinary improper-integral limit.

\subsubsection{Reduction in the $z$ variable}\label{subsubsec:vppois_z}
It remains to evaluate the $z$-integrals in the two contributions $\vpldd$ and
$\vplpv$. These are of three kinds: the elementary moment
$\int_0^{\sqrt{\veps}} z^{m+1}\,dz$ from the $\delta$-term
\eqref{eq:vppois_int6_delta}, and the two families
$\int_0^{\sqrt{\veps}} z^{n+1} I_m(r/z)\,dz$ and
$\int_0^{\sqrt{\veps}} z^{n+1}\erf(r/2z)\,dz$ from the principal value term
\eqref{eq:vppois_int8_pv}. We treat each in turn and write the result through the
scaled variable $c=r/\sqrt{\veps}$, so that every term carries an explicit factor
$\veps^{(m+2)/2}$.

The $\delta$-term involves no special functions. Carrying out
$\int_0^{\sqrt{\veps}} z^{m+1}\,dz = \veps^{(m+2)/2}/(m+2)$ in
\eqref{eq:vppois_int6_delta} gives the single sum
\begin{equation}\label{eq:vppois_delta_closed}
  \vpldd[f](\vx) = \sum_{m=0}^{N} \frac{\zeta^{\delta}_{m}}{m+2}\,\veps^{(m+2)/2}.
\end{equation}

For the principal value term, the moments $I_m$ reduce exactly as in the single
layer case \cref{subsubsec:slppois_z}: inserting the Hermite representation
\eqref{eq:I_n_exp} of $I_m$ and substituting $\omega = r/(2z)$ turns each integral
$\int_0^{\sqrt{\veps}} z^{n+1} I_m(r/z)\,dz$ into a finite combination of the
functions $Q_p$ of \eqref{eq:int_inc_gamma_def}. The error-function integrals are however new. Integrating by parts and using the same substitution gives
\begin{equation}\label{eq:vppois_erf_int}
  \int_0^{\sqrt{\veps}} z^{n+1}\erf\!\left(\frac{r}{2z}\right)dz
  = \frac{\veps^{(n+2)/2}}{n+2}\left[\,\erf\!\left(\frac{c}{2}\right)
  + \frac{2}{\sqrt{\pi}}\left(\frac{c}{2}\right)^{n+2}
    Q_{n+2}\!\left(\frac{c}{2}\right)\right],
\end{equation}
whose leading term is an error function rather than a $Q_p$, and is kept as such.
Collecting both families, the principal value part is
\begin{equation}\label{eq:vppois_int10}
\begin{split}
  \vplpv[f](\vx)
  &= \sum_{\substack{m<N\\ n<N}} \hat{\zeta}^{\pv}_{m,n}\,
     \frac{i^{n} c^{m+2}}{2^{m+n+2}}\,\veps^{(m+2)/2}
     \sum_{k=0}^{n} h_{n,k}\,Q_{m+3-k}\!\left(\frac{c}{2}\right) \\
  &\quad + \sum_{n<N} \frac{\tilde{\zeta}^{\pv}_{n}}{n+2}\,\veps^{(n+2)/2}
     \left[\,\erf\!\left(\frac{c}{2}\right)
     + \frac{2}{\sqrt{\pi}}\left(\frac{c}{2}\right)^{n+2}
       Q_{n+2}\!\left(\frac{c}{2}\right)\right].
\end{split}
\end{equation}

The local part of the volume potential is the sum of the two contributions,
\begin{equation}\label{eq:vppois_int11}
  \vpl[f](\vx) = \vpldd[f](\vx) + \vplpv[f](\vx) + \ordo{\veps}{(N+2)/2},
\end{equation}
plus an error term exponentially small in $\veps$. Collecting powers of $\veps$ yields the asymptotic
expansion of the local part of the volume potential.
\subsubsection{Obtaining $\vpl$ to a given order}\label{subsubsec:vppois_ordo}
Assume we seek asymptotic formulas for $\vpl[f](\vx)$ with an error of
$\ordo{\veps}{(P+1)/2}$ for a given $P \geq 0$. The factor of $z$ in
\eqref{eq:vppois} implies that a term $z^p$ in the integrand contributes
$\ordo{\veps}{(p+2)/2}$, so to obtain all terms up to order $\veps^{P/2}$ we
retain the integrand terms of degree $p \leq N$, with $N = P - 2$.

Collecting the closed contributions \eqref{eq:vppois_delta_closed} and
\eqref{eq:vppois_int10}, both of order $\veps^{(m+2)/2}$ with $m\ge0$, by powers of
$\veps$ gives
\begin{equation}\label{eq:vplgeneral}
  \vpl[f](\vx) = \sum_{p=0}^{P} \veps^{p/2} A_p(\vx) + \ordo{\veps}{(P+1)/2},
\end{equation}
with $A_0 = A_1 = 0$. For $P = 4$ the coefficients are
\begin{align}\label{eq:vplocO5full}
  A_2 &= \frac{1}{4} \left(2 \erf\left(\frac{c}{2}\right)
        +\frac{2 e^{-\frac{c^2}{4}} c}{\sqrt{\pi }}
        - c^2 \erfc\left(\frac{c}{2}\right) +2\right) f, \\
  A_3 &= \frac{ \left((4-2 c^2)e^{-\frac{c^2}{4}}
        +\sqrt{\pi } c^3 \, \erfc\left(\frac{c}{2}\right)\right)
        \left(\kappa_{b}  f - 2 f^{(0,1)} \right)}{12 \sqrt{\pi }}, \\
  A_4 &= \frac{1}{48}\left(\erfc\left(\frac{c}{2}\right)
        \left(4 c^4 \kappa_{b}  f^{(0,1)}-3 \left(c^4+4\right) f^{(0,2)}
        +\left(c^4-12\right) f^{(2,0)}-3 c^4 \kappa_{b} ^2 f\right)\right.\notag\\
      &\left.\quad+\frac{2 c \left(c^2-2\right) e^{-\frac{c^2}{4}}
        \left(-4 \kappa_{b}  f^{(0,1)}+3 f^{(0,2)}-f^{(2,0)}+3 \kappa_{b} ^2 f\right)}{\sqrt{\pi }}
        +24 \left(f^{(2,0)}+f^{(0,2)}\right)\right).
\end{align}
in agreement with \cite{lghtwght_ps2026}. The derivatives of $f$ are evaluated at the target point $\vx$, in the local tangent and normal coordinates at $\vb$. In the free-space case they are taken with respect to the Cartesian coordinates.

In the free-space case, the boundary is absent, corresponding to the limit $r\to\infty$ with $\veps$ fixed. This case has a separate derivation, obtained by Taylor expanding the source density at the target point. We do not include that derivation here, but note that the same coefficients are obtained by taking the limit $r\to\infty$ in each final coefficient $A_p$ in \eqref{eq:vplocO5full}.

For targets on the boundary, corresponding to $r=0$, the expansion can also be derived separately, as in the discussion following the single-layer expansion. The intermediate representation involving $Q_p(c/2)$ is introduced for $c>0$, but the preceding Fourier--Gaussian moment representation remains well defined at $r=0$, so the boundary expansion can be derived from that representation. We do not repeat that derivation here, but note that the same boundary coefficients are obtained by taking the limit $r\to 0^+$ in each final coefficient $A_p$ in \eqref{eq:vplocO5full}.

\section{The Poisson equation in $\R^3$}\label{sec:poissonr3}
We show in this section that the procedure from \cref{sec:mainresult} for the Poisson equation in two dimensions carries over naturally to the three dimensional setting. That is, we consider the Poisson equation $-\Delta u = F$ in $\Rd{3}$ with free-space Green's function $G(\vx) = \frac{1}{4\pi \abs{\vx}}$ and seek to evaluate the single layer potential, double layer potential and volume potential 
\begin{align}
  \label{eq:ps3d_single_layer}
  \slp[\sigma](\vx) &= \frac{1}{4\pi} \int_\bdry\frac{\sigma(\vxp)}{\abs{\vx-\vxp}}  \, ds(\vxp),\\
    \label{eq:ps3d_double_layer}
  \dlp[\mu](\vx) &= \frac{1}{4\pi} \int_\bdry \frac{(\vx - \vxp) \cdot \vnu(\vxp)}{\abs{\vx - \vxp}^3} \mu(\vxp)  \, ds(\vxp),\\
  \label{eq:ps3d_volume_layer}
    \vp[f](\vx) &= \frac{1}{4\pi} \int_\dom\frac{ f(\vxp)}{\abs{\vx-\vxp}}  \, d\vxp,
\end{align}
using $\ul$ and $\uh$ from \cref{eq:uL_def,eq:uH_def} for $d=3$.

The parabolic regularization strategy applies directly in the three-dimensional setting. In contrast to the two-dimensional case, there is no zero-mode ambiguity, as the Fourier representation $\fu(\vxi)=\fF(\vxi)/\abs{\vxi}^2$ is locally integrable at $\vxi=\vzero$, and the solution is uniquely determined without the need for free-space normalization. Equivalently, the three-dimensional Green’s function decays at infinity, so no logarithmic growth is present.

For the geometric setting, let $\vx$ be a target point and $\vb$ the closest point on the surface $\bdry$. We introduce local orthonormal coordinates $(x_1,x_2)$ in the tangent plane at $\vb$, aligned with the principal directions and oriented so that the normal is outward pointing from $\bdry$. In these coordinates, the surface can be represented locally as a graph over the tangent plane,
\begin{equation}\label{eq:loc_param_surface}
\vga(x_1,x_2) = (x_1,x_2,\curve(x_1,x_2)), \qquad (x_1,x_2) \in I_b,
\end{equation}
for a sufficiently small neighborhood $I_b$ of $\vb$.

In the case of layer potentials, we have $\vb = \vzero$ in the tangent plane parameterization, so $\curve(0,0) = 0$ and $\curve^{(1,0)}(0,0) = \curve^{(0,1)}(0,0) = 0$, and $\gamma^{(1,1)}(0,0)= 0$, $\gamma^{(2,0)}(0,0) = \kappa_1$ and $\curve^{(0,2)}(0,0) = \kappa_2$ where $\kappa_1$ and $\kappa_2$ are the principal curvatures of $\bdry$ at $\vb$. 
In these coordinates, $\vx = (0,0,r)$.  
For volume potentials we instead set $\vx = (0,0,0)$ with $\vb = (0,0,-r)$ and 
$\curve(0,0) = -r$.

The single layer potential is again defined as $\slpl[\sigma] := \ul[\sigma\delta(g)]$, the double layer potential as  $\dlpl[\mu] := \ul[-\nabla\cdot (\vnu\mu\delta(g))]$ and the volume potential as $\vp[f] := \ul[f\ind]$. They have the form
\begin{align}\label{eq:slppois3d_int1}
  \begin{split}
    \slpl[\sigma](\vx) &= \frac{1}{4\pi^{3}} \int_{0}^{\sqrt{\veps}} \int_{\R^3}e^{-|\veta|^2}\int_{\R^2}\Psi(\vu z) \sigma(\vga(u_1 z,u_2 z))\sqrt{1+\abs{\nabla \curve(u_1z,u_2z)}^2}\\[0.6em]
                       &e^{-i\eta_1 u_1 -i\eta_2 u_2} e^{-i\eta_3 \curve(u_1 z,u_2 z)/z}\,du_1\,du_2 \, e^{i\eta_3 r /z}\,d\veta \, dz,
  \end{split}
  \\[1em]  \begin{split}\label{eq:dlppois3d_int1}
    \dlpl[\mu](\vx) &= \frac{1}{4\pi^{3}} \int_{0}^{\sqrt{\veps}}\frac{1}{z} \int_{\R^3}i\,e^{-|\veta|^2}\veta\cdot\int_{\R^2}\Psi(\vu z)\tilde{\vnu}(\vga(u_1z,u_2z))\,\mu(\vga(u_1 z,u_2 z))e^{-i\eta_1 u_1 -i\eta_2 u_2} \\[0.6em]
                    & e^{-i\eta_3 \curve(u_1 z,u_2 z)/z}\,du_1\,du_2 \, e^{i\eta_3 r /z}\,d\veta \, dz,
  \end{split}                      
  \\[1em]  \begin{split}\label{eq:vppois3d_int1}
    \vpl[f](\vx)  &= \frac{1}{4\pi^3}\int_0^{\sqrt{\veps}} z \int_{\R^3}e^{-|\veta|^2}\int_{\R^3}\Psi(\vu z)f(\vu z)\,\theta((u_3-\curve(u_1z,u_2z)/z))e^{-i\veta\cdot\vu}e^{i\veta\cdot\vx/z}\,d\vu\,d\veta \,dz
                      \end{split}
\end{align}
with an error of $\ordo{e}{-d_B^2/(4\veps)}$. Here we have made the standard change of variables $t=z^2$, $\veta = \vxi z$ and $\vu = \vxp/z$, and with $\tilde{\vnu} = -(\curve^{(1,0)},\curve^{(0,1)}, -1)$. 

Deriving asymptotic expansions for \cref{eq:slppois3d_int1,eq:dlppois3d_int1,eq:vppois3d_int1} follows the steps outlined in \cref{sec:mainresult}. The difference lies in the presence of additional terms that must be tracked, and in the fact that the reductions in \cref{subsubsec:reduce_ueta,subsubsec:vppois_deta1_du1} must now be carried out twice, once for the pair $(u_1,\eta_1)$ and once for the pair $(u_2,\eta_2)$.

Applying this procedure yields the asymptotic expansion of the local part of the single layer potential. In terms of the mean curvature $H=\frac{\kappa_1+\kappa_2}{2}$ and the Gaussian curvature $K=\kappa_1\kappa_2$ at $\vb$, it can be written as
  \begin{equation}
    \begin{split}
      \label{eq:ps3d_single_layer_loc}
      \slpl[\sigma](\vx) &= \sqrt{\veps}\frac{\sigma_0}{2}\left(\frac{2e^{-c^2/4}}{\sqrt{\pi}} - c\erfc\left(\frac{c}{2}\right)\right) + \veps \frac{c \sigma_0 H}{2} \left(\frac{2e^{-c^2/4}}{\sqrt{\pi}} - c\erfc\left(\frac{c}{2}\right)\right) \\
      &+ \frac{\veps^{3/2}}{48\sqrt{\pi}}
\Bigl[
\bigl((16+16c^2)e^{-c^2/4}-8c^3\sqrt{\pi}\,\erfc(c/2)\bigr)H^2\,\sigma_0\\
&-\bigl((4+4c^2)e^{-c^2/4}-2c^3\sqrt{\pi}\,\erfc(c/2)\bigr)K\,\sigma_0+\bigl((16-8c^2)e^{-c^2/4}+4c^3\sqrt{\pi}\,\erfc(c/2)\bigr)\Delta_\Gamma \sigma_0
  \Bigr]\\
      & + \ordo{\veps}{2},
    \end{split}
  \end{equation}    
where $c = r/\sqrt{\veps}$ and $\Delta_\bdry$ denotes the Laplace--Beltrami operator on $\bdry$. In local orthonormal coordinates at the base point $\vb$ it satisfies
\begin{equation}
\Delta_\bdry \sigma_{b} = \partial_{x_1}^{2}\sigma(\vb) + \partial_{x_2}^{2}\sigma(\vb).
\end{equation}

For the double layer potential we have
\begin{equation}
  \begin{split}
    \dlpl[\mu](\vx) &= -\frac{1}{2} \erfc\left(\frac{c}{2}\right) \mu_b - \sqrt{\veps}\,\frac{e^{-c^2/4}}{\sqrt{\pi}}\,H\,\mu_b \\
                    &+ \veps\left(
                      \frac{c^2}{4}\,\erfc\!\left(\frac{c}{2}\right)\Delta_\bdry \mu_b
                      -\frac{ce^{-c^2/4}}{2\sqrt{\pi}}
                      \Bigl(
                      (3H^2-K)\mu_b+\Delta_\bdry \mu_b
                      \Bigr)
                      \right) \\
                    &-\veps^{3/2}\Biggl[
-\frac{e^{-c^2/4}}{2\sqrt{\pi}}
\Bigl(
(c^2-2)H(5H^2-3K)\mu_b+\mu_b\,\Delta_\bdry H
\Bigr)\\
&+\left(
\frac{c^3}{12}\,\erfc\!\left(\frac{c}{2}\right)
-\frac{(1+c^2)e^{-c^2/4}}{6\sqrt{\pi}}
\right)\Bigl(
(\kappa_1+3\kappa_2)\mu_b^{(0,2)}+(3\kappa_1+\kappa_2)\mu_b^{(2,0)}
\Bigr)\\
&+\left(
\frac{c^3}{12}\,\erfc\!\left(\frac{c}{2}\right)
-\frac{(4+c^2)e^{-c^2/4}}{6\sqrt{\pi}}
\right)
\Bigl(
\mu_b^{(0,1)}(\gamma_b^{(0,3)}+\gamma_b^{(2,1)})+\mu_b^{(1,0)}(\gamma_b^{(1,2)}+\gamma_b^{(3,0)})
\Bigr)\Biggr]\\
                    &+ \ordo{\veps}{2}.
  \end{split}
\end{equation}
To order $\veps$, the local asymptotic correction can be expressed entirely in terms of invariant surface quantities: the Laplace--Beltrami operator applied to the density, the mean curvature, and the Gaussian curvature. At higher orders, the dependence on the local parameterization does not disappear, but reorganizes into invariant geometric quantities. 

Finally, we have the volume potential,
\begin{equation}
  \begin{split}
    \vpl[f](\vx) &= \frac{\veps}{4}\left(2\erf\left(\frac{c}{2}\right)+\frac{2 e^{-\frac{c^2}{4}}c}{\sqrt{\pi}}-c^2\erfc\left(\frac{c}{2}\right)+2\right)f\\
                 &- \frac{\veps^{3/2}}{12\sqrt{\pi}}\left((4-2c^2)e^{-\frac{c^2}{4}} + \sqrt{\pi}c^3\erfc\left(\frac{c}{2}\right)\right)(2Hf-2f^{(0,0,1)})\\
                 &+\veps^{2}\Biggl(\frac{1}{12} c(c^2 - 2) e^{-c^2/4} \sqrt{\pi}\Bigl((12H^2-4K) f- 8H f^{(0,0,1)} + 3 f^{(0,0,2)} - f^{(0,2,0)} - f^{(2,0,0)}\Bigr) \\
    &+ \pi \Delta f + \frac{\pi}{24} \erfc\left(\frac{c}{2}\right)
    \Bigl( c^4 \bigl(
            -(12H^2- 4K) f
            + 8H f^{(0,0,1)}
    \\
    &- 3 f^{(0,0,2)}
            + f^{(0,2,0)}
            + f^{(2,0,0)}
        \bigr)
        - 12\Delta f
    \Bigr)
      \Biggr)\\
    &+ \ordo{\veps}{5/2}.
    \end{split}
\end{equation}
We see that each term in the two-dimensional expansion \eqref{eq:vplocO5full} has a direct three-dimensional analog obtained by replacing scalar curvature by the principal curvatures, tangential derivatives by surface derivatives, and higher derivatives of the local parameterization by derivatives of the curvature.

\section{Coupled Strongly Elliptic Systems}\label{sec:coupsys}
We extend the parabolic regularization strategy from \cref{sec:regstrat} to coupled strongly elliptic systems. As for a single equation, the construction depends only on the Fourier symbol matrix $\fM(\vxi)$ of the (matrix) differential operator $M$, and does not require a Green's function in physical space. 
In this paper, we will assume that
the eigenvalues of $\fM(\vxi)$ have strictly positive real part for $\vxi \neq \vzero$. 
In that case, the matrix exponential
$e^{-\fM(\vxi)t}$ decays as $t \to \infty$ and the localization mechanism of \cref{sec:regstrat} carries over directly. Thus the geometric setup, localization arguments, and sequence of reductions are identical to those of the preceding sections. 
The only new ingredient is a more involved analysis of the matrix exponential 
to account for the coupling between system components.

As a model problem, we consider the system of constant-coefficient PDEs
\begin{equation}\label{eq:coupsys}
  \begin{cases}
    -\Delta v_1 +  \va\cdot\nabla v_2 + a_3 v_2  &= F_1\\
    -\Delta v_2 +  \vb\cdot\nabla v_1  + b_3 v_1   &= F_2
  \end{cases}
  \quad \text{ in } \R^2,
\end{equation}
where real valued $\va=(a_1,a_2)$ and $\vb=(b_1,b_2)$ introduce first-order coupling, and $a_3,b_3\in\R$ introduce zeroth-order coupling.

Let $\vv = (v_1,v_2)$ and $\vF = (F_1,F_2)$, so that we can rewrite \cref{eq:coupsys} in the form
\begin{equation}
  M \vv = \vF \quad \text{in }\R^2,
\end{equation}
where the matrix-valued differential operator is
\begin{equation}\label{eq:diff_mat}
  M := \begin{bmatrix}
   -\Delta & a_1\partial_{x_1} + a_2\partial_{x_2} + a_3\\
    b_1\partial_{x_1} + b_2\partial_{x_2}  + b_3& -\Delta
  \end{bmatrix}.
\end{equation}
We define the Fourier symbol-matrix $\fM(\vxi)$ associated with the differential operator $M$ by replacing each differential operator componentwise by its Fourier symbol, that is, $\partial_{x_j} \mapsto i\xi_j$ and $-\Delta \mapsto \abs{\vxi}^2$, yielding
\begin{equation}\label{eq:symb_mat}
  \fM(\vxi) := \begin{bmatrix}
    \abs{\vxi}^2 & i(\va\cdot\vxi) + a_3\\
    i(\vb\cdot\vxi) + b_3 & \abs{\vxi}^2
  \end{bmatrix}.
\end{equation}

The system \eqref{eq:coupsys} is strongly elliptic, since its principal symbol $\abs{\vxi}^2 I$ is positive definite for all $\vxi \neq \vzero$, regardless of the lower-order parameters $\va,a_3,\vb,b_3$. The eigenvalues of $\fM(\vxi)$ are
\begin{equation}
    \lambda_\pm(\vxi) = \abs{\vxi}^2 \pm \sqrt{(i\va\cdot\vxi + a_3)(i\vb\cdot\vxi + b_3)}.
\end{equation}

We assume that the coefficients $\va,a_3,\vb,b_3$ are chosen so that $\fM(\vxi)$ is invertible for all $\vxi\in\R^2$, and that its eigenvalues $\lambda_\pm(\vxi)$ satisfy $\real(\lambda_\pm(\xi))>0$, for nonzero $\vxi$. The positivity of the real parts ensures that $e^{-\fM(\vxi)t}$ decays as $t\to\infty$ for each $\vxi\neq \vzero$, and allows the localization mechanism of \cref{sec:regstrat} to carry over directly. Since $\fM(\vzero)$ is
invertible, $\fM(\vxi)^{-1}$ is smooth at the origin and no zero-mode correction is needed, in contrast to \cref{sec:regstrat}. For the present system \cref{eq:symb_mat}, the condition $a_3 b_3<0$ is sufficient for invertibility of $\fM(\vzero)$. In the special case $\va b_3+\vb a_3=\vzero$ and $a_3 b_3<0$, the product  $(i\va\cdot\vxi+a_3)(i\vb\cdot\vxi+b_3)$ is negative for all $\vxi\in\R^2$, so that $\real(\lambda_\pm(\vxi)=\abs{\vxi}^2>0$ for all $\vxi\neq \vzero$, and the conditions these are automatically satisfied.

Taking the Fourier transform of \cref{eq:coupsys}, we obtain $\fM(\vxi)\fhat{\vv}(\vxi) = \fhat{\vF}(\vxi)$ for each $\vxi\in\R^2$. The construction of \cref{sec:regstrat} extends directly to this setting, yielding the vector-valued potential $\vpot[\vF] := \vpoth[\vF] + 
\vpotl[\vF]$, with history part
\begin{equation}\label{eq:vpothdef}
\vpoth[\vF](\vx) := \frac{1}{(2\pi)^{2}}\int_{\R^2}\int_{\veps}^{\infty}
e^{-\fM(\vxi)t}\int_{\R^2}\vF(\vxp)e^{-i\vxi\cdot\vxp}\,d\vxp\,dt\,e^{i\vxi\cdot\vx} d\vxi,
\end{equation}
and local part
\begin{equation}\label{eq:vpotldef}
  \vpotl[\vF](\vx) := \frac{1}{(2\pi)^{2}}\int_0^{\veps}\int_{\R^2}
  e^{-\fM(\vxi)t}\int_{\R^2}\vF(\vxp)e^{-i\vxi\cdot\vxp}\,d\vxp e^{i\vxi\cdot\vx}\,d\vxi\,dt.  
\end{equation}
Whenever $\fM(\vxi)$ is invertible, the local part admits the closed form
\begin{equation}
\int_0^{\veps} e^{-\fM(\vxi)t}\,dt = \fM(\vxi)^{-1}\bigl(I - e^{-\fM(\vxi)\veps}\bigr).
\end{equation}
By assumption, $\fM(\vxi)$ is invertible for all $\vxi \in \R^2$, and its eigenvalues $\lambda_\pm(\vxi)$ satisfy $\real(\lambda_\pm(\vxi)) > 0$ for all $\vxi \neq \vzero$. The positivity of the real parts implies that the matrix exponential $e^{-\fM(\vxi)\veps}$ decays exponentially as $\abs{\vxi}\to\infty$. Moreover, since $\fM(\vxi)^{-1}$ is smooth away from $\vxi=\vzero$ and the Fourier transform of the data is smooth, the apparent singularity at $\vxi=\vzero$ is removable.

The asymptotic analysis proceeds as in \cref{sec:mainresult}. After the change 
of variables $z = \sqrt{t}$ and $\veta = z\vxi$, the Fourier symbol $\fM(\vxi)$ separates into its leading isotropic part and lower-order coupling terms. In particular, the Laplacian contributes the dominant term $\abs{\vxi}^2 I$, while the first- and zeroth-order coupling terms are suppressed by factors of $z$ and $z^2$, respectively. This yields the factorization
\begin{equation}\label{eq:factorexpfM}
    e^{-\fM(\veta/z)z^2} = e^{-\abs{\veta}^2 I} e^{-R(\veta,z)},
\end{equation}
where the remainder $R(\veta,z)$ collects all lower-order contributions and is given by
\begin{equation}
  R(\veta,z) = iz\begin{bmatrix}
    0 & (\va \cdot \veta)\\
    (\vb \cdot \veta) & 0
  \end{bmatrix}
  + z^2\begin{bmatrix}
    0 & a_3\\
    b_3 & 0
  \end{bmatrix}.
\end{equation}
The matrix $R(\veta,z)$ is $O(z)$ as $z \to 0^+$, with a linear contribution in $\veta$ arising from the first-order coupling and a constant contribution from the zeroth-order terms. Expanding $e^{-R(\veta,z)}$ in powers of $z$ therefore produces a systematic expansion in which the components are coupled order by order, with each power of $z$ introducing higher-order interactions between $v_1$ and $v_2$.

We now derive asymptotic expansions for the local parts of the three potential types. The volume potential is treated first in \cref{subsec:vpsys}, as the two-dimensional source geometry introduces the most structure. The single and double layer potentials are considered in \cref{subsec:slpsys,subsec:dlpsys}, where the reduction follows the same pattern with the Heaviside factor replaced by a layer density on $\bdry$.

\subsection{The Volume Potential $\vp$}\label{subsec:vpsys}
We define the vector-valued volume potential $\bm{\calV}[\vf] = [\calV^{(1)}[\vf] \;\, \calV^{(2)}[\vf]] := \vpot[\vf\ind]$, where $\vf = [f_1 \; f_2]$ with source densities $f_1,f_2$ smooth on $\dom$, and $\ind$ is the indicator function of $\dom$. The local part $\bm{\calV}_L[\vf]$ and history part $\bm{\calV}_H[\vf]$ are defined analogously through \cref{eq:vpotldef,eq:vpothdef}.

The history part is given by
\begin{equation}
\bm{\calV}_H[\vf](\vx) := \frac{1}{(2\pi)^{2}}\int_{\R^2}\fM^{-1}(\vxi)e^{-\fM(\vxi)\veps}\int_{\R^2}\vf(\vxp)\ind(\vxp)e^{-i\vxi\cdot\vxp}\,d\vxp e^{i\vxi\cdot\vx}\,d\vxi,
\end{equation}
and requires no zero-mode correction since $\fM(0)$ is invertible. We make no further comments on computing the history part for the case of coupled PDEs. We simply observe that the matrix-valued multiplier $ \fhat{M}(\vxi)^{-1}e^{-\veps \fhat{M}(\vxi)}$ is applied mode by mode in Fourier space, after which the inverse transform is evaluated componentwise using scalar PDE techniques \cite{vico,dmk2024}.

Under the same assumptions and geometric localization as in \cref{subsec:vppois} we write \cref{eq:vpotldef} as
\begin{equation}\label{eq:vppois_sys}
  \bm{\calV}_L[\vf](\vx) = \frac{1}{2\pi^2}\int_0^{\sqrt{\veps}}z\int_{\R^2}e^{-\abs{\veta}^2 I}e^{-R(\veta,z)}\int_{\R^2}\Psi(\vu z)\vf(\vu z)\,\theta\!
  \left(u_2-\curvex(u_1,z)\right)e^{-i\veta\cdot\vu}\,d\vu\,d\veta\,dz,
\end{equation}
with an error of $\ordo{e}{-d_B^2/(4\veps)}$, where $\vu z = \vxp$ and $\curvex(u_1,z)$ from \eqref{eq:curve_tail}. We proceed as in \cref{subsec:vppois} by expanding $\curvex$ and $\vf$ in Taylor series,
\begin{equation}
    \vf(\vu z)
  = \sum_{m+n\le N} \frac{1}{m!\,n!} u_1^m u_2^n z^{m+n}\vf^{(m,n)}  + R_{N+1}(\vu z),
\end{equation}
where $\vf^{(m,n)} = [\partial_1^{\,m}\partial_2^{\,n} f_{1}\; \partial_1^{\,m}\partial_2^{\,n} f_{2}]$ and $f_{1} = f_{1}(\vx)$, $f_{2} = f_{2}(\vx)$. The matrix-valued exponential $e^{-R(\veta,z)}$ is expanded as
\begin{equation}
e^{-R(\veta,z)} = \sum_{n=0}^{N}\frac{(-1)^n}{n!}(R(\veta,z))^n +R_{N+1}(\veta z),
\end{equation}
yielding finite sums of monomials in $\eta_1$ and $\eta_2$ to order $N$. The reductions then proceed exactly as in \cref{subsec:vppois}, with the expansion of $e^{-R(\veta,z)}$ providing the coupling between components. Straightforward calculations, similar to those in \eqref{subsubsec:reduce_ueta}, shows that $R_{N+1}$ after integration in $z$ introduces an error of order $\ordo{\veps}{(N+2)/2}$.

Let $P \geq 0$ be given. To obtain an expansion of $\bm{\calV}_L[\vf](\vx)$ with error $\ordo{\veps}{(P+1)/2}$, it suffices to expand each Taylor series to order $N = P-2$, exactly as in \cref{subsec:vppois}. Since $R(\veta,z)$ contains only positive powers of $z$, the same power counting as for $\vpldd$ in \cref{subsec:vppois} applies here as well. It follows that the asymptotic expansion takes the form
\begin{equation}
\bm{\calV}_L[\vf](\vx) = \sum_{p=0}^{P}\bm{A}_p\veps^{p/2} + \ordo{\veps}{(P+1)/2},
\end{equation}
with $\bm{A}_p = [A^{(1)}_p\;A^{(2)}_p]$.

For $P = 4$ we present the coefficients $\bm{A}_p$. The extra factor of $z$ in \cref{eq:vppois_sys} means the lowest-order contribution to the integrand is $O(z)$, which after integration yields $O(\veps)$. Consequently $\bm{A}_0 = \bm{0}$ and $\bm{A}_1 = \bm{0}$, as in \cref{subsec:vppois}. Let $c = r/\sqrt{\veps}$, and note that in the local coordinate system used below, $x_1$ is the tangential coordinate and $x_2$ is the inward normal coordinate. Thus, if $\vnu$ denotes the outward unit normal and $\bm{\tau}$ the unit tangent, then
\begin{equation}
\partial_{x_1}=\partial_{\bm{\tau}} ,
\qquad
\partial_{x_2}=-\partial_{\vnu} .
\end{equation}
Consequently, the coefficients appearing in the local formulas are
\begin{equation}
\alpha_1=\va\cdot\bm{\tau},
\qquad
\alpha_2=-\va\cdot\vnu,
\qquad
\alpha_3=a_3,
\qquad
\beta_1=\vb\cdot\bm{\tau},
\qquad
\beta_2=-\vb\cdot\vnu.
\qquad
\beta_3=b_3.
\end{equation}
We have
\begin{align}
  A^{(1)}_2
  &= \frac{1}{4}\left(
      4 + \frac{2c e^{-c^2/4}}{\sqrt{\pi}}
      - (2+c^2)\erfc\left(\frac{c}{2}\right)
    \right) f_1
\\
  A^{(1)}_3
  &= -\frac{\alpha_2}{12}\left(
      \frac{2(c^2-2)e^{-c^2/4}}{\sqrt{\pi}}
      - c^3\erfc\left(\frac{c}{2}\right)
    \right) f_2\\
  &- \frac{e^{-c^2/4}}{12\sqrt{\pi}}
    \left(
      4 - 2c^2 + c^3 e^{c^2/4}\sqrt{\pi}\erfc\left(\frac{c}{2}\right)
    \right)
    \left(\kappa_b f_1 - 2 f_1^{(0,1)}\right)
\\
  A^{(1)}_4
  &= \frac{1}{48}\Biggl[
      \frac{2c(c^2-2)e^{-c^2/4}}{\sqrt{\pi}}
      \left(
        3\kappa_b^2 f_1 - 4\kappa_b f_1^{(0,1)} + 3f_1^{(0,2)} - f_1^{(2,0)}
      \right)
      + 24\left(f_1^{(0,2)} + f_1^{(2,0)}\right)
\notag\\
&\qquad
      - \erfc\left(\frac{c}{2}\right)
      \left[
        c^4
        \left(
          3\kappa_b^2 f_1 - 4\kappa_b f_1^{(0,1)} + 3f_1^{(0,2)} - f_1^{(2,0)}
        \right)
        + 12\left(f_1^{(0,2)} + f_1^{(2,0)}\right)
      \right]
\notag\\
&\qquad
      + \alpha_2 \beta_2
      \left(
        \frac{2c(c^2-2)e^{-c^2/4}}{\sqrt{\pi}}
        - c^4 \erfc\left(\frac{c}{2}\right)
      \right) f_1
\notag\\
&\qquad
      + \alpha_3
      \left(
        -24
        - \frac{2c(2-c^2)e^{-c^2/4}}{\sqrt{\pi}}
        + (12-c^4)\erfc\left(\frac{c}{2}\right)
      \right) f_2
\notag\\
&\qquad
      + \alpha_2 \kappa_b
      \left(
        c^4 \erfc\left(\frac{c}{2}\right)
        - \frac{2c(c^2-2)e^{-c^2/4}}{\sqrt{\pi}}
      \right) f_2
\notag\\
&\qquad
      + \alpha_2
      \left(
        \frac{3c(c^2-2)e^{-c^2/4}}{\sqrt{\pi}}
        + 3\left(
          -8 + (4+c^4)\erfc\left(\frac{c}{2}\right)
        \right)
      \right) f_2^{(0,1)}
\notag\\
&\qquad
      + \alpha_1
      \left(
        -24 + (12-c^4)\erfc\left(\frac{c}{2}\right)
        + \frac{2c-c^3}{\sqrt{\pi}}e^{-c^2/4}
      \right) f_2^{(1,0)}
    \Biggr].
\end{align}
and
\begin{align}
  A^{(2)}_2
  &= \frac{1}{4}\left(
      4 + \frac{2c e^{-c^2/4}}{\sqrt{\pi}}
      - (2+c^2)\erfc\left(\frac{c}{2}\right)
    \right) f_2
\\
  A^{(2)}_3
  &= -\frac{\beta_2}{12}\left(
      \frac{2(c^2-2)e^{-c^2/4}}{\sqrt{\pi}}
      - c^3\erfc\left(\frac{c}{2}\right)
    \right) f_1\\
  &- \frac{e^{-c^2/4}}{12\sqrt{\pi}}
    \left(
      4 - 2c^2 + c^3 e^{c^2/4}\sqrt{\pi}\erfc\left(\frac{c}{2}\right)
    \right)
    \left(\kappa_b f_2 - 2 f_2^{(0,1)}\right)
\end{align}
\begin{align}
  A^{(2)}_4
  &= \frac{1}{48}\Biggl[
      \frac{2c(c^2-2)e^{-c^2/4}}{\sqrt{\pi}}
      \left(
        3\kappa_b^2 f_2 - 4\kappa_b f_2^{(0,1)} + 3f_2^{(0,2)} - f_2^{(2,0)}
      \right)
      + 24\left(f_2^{(0,2)} + f_2^{(2,0)}\right)
\notag\\
&\qquad
      - \erfc\left(\frac{c}{2}\right)
      \left[
        c^4
        \left(
          3\kappa_b^2 f_2 - 4\kappa_b f_2^{(0,1)} + 3f_2^{(0,2)} - f_2^{(2,0)}
        \right)
        + 12\left(f_2^{(0,2)} + f_2^{(2,0)}\right)
      \right]
\notag\\
&\qquad
      + \alpha_2 \beta_2
      \left(
        \frac{2c(c^2-2)e^{-c^2/4}}{\sqrt{\pi}}
        - c^4 \erfc\left(\frac{c}{2}\right)
      \right) f_2
\notag\\
&\qquad
      + \beta_3
      \left(
        -24
        - \frac{2c(2-c^2)e^{-c^2/4}}{\sqrt{\pi}}
        + (12-c^4)\erfc\left(\frac{c}{2}\right)
      \right) f_1
\notag\\
&\qquad
      + \beta_2 \kappa_b
      \left(
        c^4 \erfc\left(\frac{c}{2}\right)
        - \frac{2c(c^2-2)e^{-c^2/4}}{\sqrt{\pi}}
      \right) f_1
\notag\\
&\qquad
      + \beta_2
      \left(
        \frac{3c(c^2-2)e^{-c^2/4}}{\sqrt{\pi}}
        + 3\left(
          -8 + (4+c^4)\erfc\left(\frac{c}{2}\right)
        \right)
      \right) f_1^{(0,1)}
\notag\\
&\qquad
      + \beta_1
      \left(
        -24 + (12-c^4)\erfc\left(\frac{c}{2}\right)
        + \frac{2c-c^3}{\sqrt{\pi}}e^{-c^2/4}
      \right) f_1^{(1,0)}
    \Biggr].
\end{align}

Several structural properties of the expansion are worth noting. The two components $A^{(1)}_p$ and $A^{(2)}_p$ are not related by a simple symmetry under exchange of $f_1$ and $f_2$, due to the asymmetric off-diagonal structure of $\fM(\vxi)$. Consequently, the expansion of $\bm{\mathcal{V}}_L^{(2)}$ cannot be obtained from that of $\bm{\mathcal{V}}_L^{(1)}$ by a direct exchange of the data, but instead requires exchanging the coefficients $(\alpha_1,\alpha_2,\alpha_3)$ with $(\beta_1,\beta_2,\beta_3)$.

At order $\veps^{3/2}$, the coupling between the two components is governed exclusively by the normal convection coefficients $\alpha_2$ and $\beta_2$. The contribution to $\bm{\mathcal{V}}_L^{(1)}$ depends only on $\alpha_2$ and involves $f_2$, while the corresponding term in $\bm{\mathcal{V}}_L^{(2)}$ depends only on $\beta_2$ and involves $f_1$. 

At order $\veps^2$, the coupling decomposes into three distinct contributions. The first arises from the constant off-diagonal terms $\alpha_3$ and $\beta_3$, which produce zeroth-order coupling. The second involves the normal convection coefficients $\alpha_2$ and $\beta_2$, which continue to couple the components through normal derivatives. The third involves the tangential convection coefficients $\alpha_1$ and $\beta_1$, which introduce coupling through tangential derivatives.

The uncoupled part at order $\veps^2$, consisting of terms involving only $f_1$ in $\mathcal{V}_L^{(1)}$ and only $f_2$ in $\mathcal{V}_L^{(2)}$, coincides in structure with the scalar volume potential expansion \eqref{eq:vplocO5full}, with $f$ replaced by $f_1$ and $f_2$, respectively. In particular, it depends only on geometric quantities such as curvature and on second derivatives of the data, and is independent of the coupling parameters. The same structure persists at lower orders, confirming that the system expansion reduces to the scalar case in the absence of off-diagonal terms.

\subsection{The Single Layer Potential $\slp$}\label{subsec:slpsys}
The vector-valued single layer potential $\bm{\calS}[\vsigma]$ is defined 
analogously to the scalar case \cref{subsec:slppois}, with density 
$\vsigma = [\sigma_1\;\sigma_2]$ supported on $\bdry$. The local part 
$\bm{\calS}_L[\vsigma]$ is obtained from \cref{eq:vpotldef} with 
$\vF = \vsigma\delta(g)$, and after the same change of variables and 
localization as in \cref{subsec:slppois}, takes the form
\begin{equation}\label{eq:slpsys}
  \bm{\calS}_L[\vsigma](\vx) = \frac{1}{2\pi^2}\int_0^{\sqrt{\veps}}
  \int_{\R^2}e^{-\abs{\veta}^2 I}e^{-R(\veta,z)}\int_{\R}\psi(uz)
  \vsigma(\curve(u z))\sqrt{1+(\curve'(u z))^2}\,
  e^{-i\eta_1 u}e^{-i\eta_2\curvex(u,z)}\,du\,
  e^{i\eta_2 r/z}\,d\veta\,dz.
\end{equation}
with an error of $\ordo{e}{-d_B^2/(4\veps)}$. The reduction proceeds exactly as in \cref{subsec:slppois}, with the 
expansion of $e^{-R(\veta,z)}$ providing the coupling between components 
order by order in $z$. The asymptotic expansion takes the form
\begin{equation}
\bm{\calS}_L[\vsigma](\vx) = \sum_{p=1}^{P}\bm{A}_p\veps^{p/2} + 
\ordo{\veps}{(P+1)/2},
\end{equation}
with $\bm{A}_p = [A^{(1)}_p\;A^{(2)}_p]$. For $P = 4$ we have
\begin{align}
A^{(1)}_{1}
&=
\sigma_{1}\left(
\frac{e^{-c^{2}/4}}{\sqrt{\pi}}-\frac c2\,\erfc\left(\frac{c}{2}\right)
\right)\\
A^{(1)}_{2}
&=
\frac c2\left(\frac{e^{-c^{2}/4}}{\sqrt{\pi}}-\frac c2\,\erfc\left(\frac{c}{2}\right)\right)\bigl(k\,\sigma_{1}+\alpha_{2}\sigma_{2}\bigr)\\
A^{(1)}_{3}
&=
\sigma_{1}\Bigg[
\frac{e^{-c^{2}/4}}{\sqrt{\pi}}\left(
-\frac{\alpha_{2}\beta_{2}}{12}
+\frac{\alpha_{2}\beta_{2}}{6}c^{2}
+\frac{k^{2}}{12}
+\frac{k^{2}}{3}c^{2}
\right)
-\erfc\left(\frac{c}{2}\right)\left(
\frac{\alpha_{2}\beta_{2}}{12}c^{3}
+\frac{k^{2}}{6}c^{3}
\right)
\Bigg]\notag\\
&\quad+
\sigma_{2}\Bigg[
\frac{e^{-c^{2}/4}}{\sqrt{\pi}}\left(
-\frac{\alpha_{3}}{3}
+\frac{\alpha_{3}}{6}c^{2}
-\frac{\alpha_{2}k}{6}
+\frac{\alpha_{2}k}{3}c^{2}
\right)
-\erfc\left(\frac{c}{2}\right)\left(
\frac{\alpha_{3}}{12}c^{3}
+\frac{\alpha_{2}k}{6}c^{3}
\right)
\Bigg]\notag\\
&\quad+
\sigma_{2}'\Bigg[
\frac{e^{-c^{2}/4}}{\sqrt{\pi}}\left(
-\frac{\alpha_{1}}{3}
+\frac{\alpha_{1}}{6}c^{2}
\right)
-\erfc\left(\frac{c}{2}\right)\frac{\alpha_{1}}{12}c^{3}
\Bigg]\notag\\
&\quad+
\sigma_{1}''\Bigg[
\frac{e^{-c^{2}/4}}{\sqrt{\pi}}\left(
\frac13-\frac{c^{2}}6
\right)
+\erfc\left(\frac{c}{2}\right)\,\frac{c^{3}}{12}
\Bigg]\\
A^{(1)}_{4}
&=
\sigma_{1}\Bigg[
\frac{e^{-c^{2}/4}}{\sqrt{\pi}}\Bigg(
-\frac{\alpha_{3}\beta_{2}}{12}c
-\frac{\alpha_{2}\beta_{3}}{12}c
+\frac{\alpha_{3}\beta_{2}}{24}c^{3}
+\frac{\alpha_{2}\beta_{3}}{24}c^{3}
+\frac{\alpha_{1}\beta_{1}k}{12}c
-\frac{\alpha_{2}\beta_{2}k}{8}c
-\frac{\alpha_{1}\beta_{1}k}{24}c^{3}\notag\\
&\qquad
+\frac{\alpha_{2}\beta_{2}k}{8}c^{3}
-\frac{k^{3}}{8}c
+\frac{3k^{3}}{8}c^{3}
\Bigg)
+\erfc\left(\frac{c}{2}\right)\Bigg(
-\frac{\alpha_{3}\beta_{2}}{48}c^{4}
-\frac{\alpha_{2}\beta_{3}}{48}c^{4}
+\frac{\alpha_{1}\beta_{1}k}{48}c^{4}
-\frac{\alpha_{2}\beta_{2}k}{16}c^{4}
-\frac{3k^{3}}{16}c^{4}
\Bigg)
\Bigg]\notag\\
&\quad+
\sigma_{2}\Bigg[
\frac{e^{-c^{2}/4}}{\sqrt{\pi}}\Bigg(
-\frac{\alpha_{2}^{2}\beta_{2}}{24}c
+\frac{\alpha_{2}^{2}\beta_{2}}{24}c^{3}
-\frac{\alpha_{3}k}{6}c
+\frac{\alpha_{3}k}{12}c^{3}
-\frac{5\alpha_{2}k^{2}}{24}c
+\frac{7\alpha_{2}k^{2}}{24}c^{3}
\Bigg)\notag\\
&\qquad
+\erfc\left(\frac{c}{2}\right)\Bigg(
-\frac{\alpha_{2}^{2}\beta_{2}}{48}c^{4}
-\frac{\alpha_{3}k}{24}c^{4}
-\frac{7\alpha_{2}k^{2}}{48}c^{4}
\Bigg)
\Bigg]\notag\\
&\quad+
\sigma_{1}'\Bigg[
\frac{e^{-c^{2}/4}}{\sqrt{\pi}}\left(
-\frac{\beta_{1}\alpha_{2}}{12}c
-\frac{\alpha_{1}\beta_{2}}{12}c
+\frac{\beta_{1}\alpha_{2}}{24}c^{3}
+\frac{\alpha_{1}\beta_{2}}{24}c^{3}
\right)
-\erfc\left(\frac{c}{2}\right)\left(
\frac{\beta_{1}\alpha_{2}}{48}c^{4}
+\frac{\alpha_{1}\beta_{2}}{48}c^{4}
\right)
\Bigg]\notag\\
&\quad+
\sigma_{2}'\Bigg[
\frac{e^{-c^{2}/4}}{\sqrt{\pi}}\left(
-\frac{\alpha_{1}k}{2}c
+\frac{\alpha_{1}k}{4}c^{3}
\right)
-\erfc\left(\frac{c}{2}\right)\frac{\alpha_{1}k}{8}c^{4}
\Bigg]\notag\\
&\quad+
\sigma_{1}''\Bigg[
\frac{e^{-c^{2}/4}}{\sqrt{\pi}}\left(
\frac{k}{2}c-\frac{k}{4}c^{3}
\right)
+\erfc\left(\frac{c}{2}\right)\,\frac{k}{8}c^{4}
\Bigg]\notag\\
&\quad+
\sigma_{2}''\Bigg[
\frac{e^{-c^{2}/4}}{\sqrt{\pi}}\left(
\frac{\alpha_{2}}{6}c-\frac{\alpha_{2}}{12}c^{3}
\right)
+\erfc\left(\frac{c}{2}\right)\,\frac{\alpha_{2}}{24}c^{4}
\Bigg]\notag\\
&\quad+
\sigma_{2}\,\gamma^{(3)}\Bigg[
\frac{e^{-c^{2}/4}}{\sqrt{\pi}}\left(
-\frac{\alpha_{1}}{6}c+\frac{\alpha_{1}}{12}c^{3}
\right)
-\erfc\left(\frac{c}{2}\right)\,\frac{\alpha_{1}}{24}c^{4}
\Bigg]\notag\\
&\quad+
\sigma_{1}'\,\gamma^{(3)}\Bigg[
\frac{e^{-c^{2}/4}}{\sqrt{\pi}}\left(
\frac13 c-\frac16 c^{3}
\right)
+\erfc\left(\frac{c}{2}\right)\,\frac1{12}c^{4}
\Bigg]\notag\\
&\quad+
\sigma_{1}\,\gamma^{(4)}\Bigg[
\frac{e^{-c^{2}/4}}{\sqrt{\pi}}\left(
\frac1{12}c-\frac1{24}c^{3}
\right)
+\erfc\left(\frac{c}{2}\right)\,\frac1{48}c^{4}
\Bigg]
\end{align}
and
\begin{align}
A^{(2)}_{1}
&=
\sigma_{2}\left(
\frac{e^{-c^{2}/4}}{\sqrt{\pi}}-\frac c2\,\erfc\left(\frac{c}{2}\right)
\right)\\
A^{(2)}_{2}
&=
\frac c2\left(\frac{e^{-c^{2}/4}}{\sqrt{\pi}}-\frac c2\,\erfc\left(\frac{c}{2}\right)\right)\bigl(\beta_{2}\,\sigma_{1}+k\,\sigma_{2}\bigr)\\
A^{(2)}_{3}
&=
\sigma_{1}\Bigg[
\frac{e^{-c^{2}/4}}{\sqrt{\pi}}\left(
-\frac{\beta_{3}}{3}
+\frac{\beta_{3}}{6}c^{2}
-\frac{\beta_{2}k}{6}
+\frac{\beta_{2}k}{3}c^{2}
\right)
-\erfc\left(\frac{c}{2}\right)\left(
\frac{\beta_{3}}{12}c^{3}
+\frac{\beta_{2}k}{6}c^{3}
\right)
\Bigg]\notag\\
&\quad+
\sigma_{2}\Bigg[
\frac{e^{-c^{2}/4}}{\sqrt{\pi}}\left(
-\frac{\alpha_{2}\beta_{2}}{12}
+\frac{\alpha_{2}\beta_{2}}{6}c^{2}
+\frac{k^{2}}{12}
+\frac{k^{2}}{3}c^{2}
\right)
-\erfc\left(\frac{c}{2}\right)\left(
\frac{\alpha_{2}\beta_{2}}{12}c^{3}
+\frac{k^{2}}{6}c^{3}
\right)
\Bigg]\notag\\
&\quad+
\sigma_{1}'\Bigg[
\frac{e^{-c^{2}/4}}{\sqrt{\pi}}\left(
-\frac{\beta_{1}}{3}
+\frac{\beta_{1}}{6}c^{2}
\right)
-\erfc\left(\frac{c}{2}\right)\frac{\beta_{1}}{12}c^{3}
\Bigg]\notag\\
&\quad+
\sigma_{2}''\Bigg[
\frac{e^{-c^{2}/4}}{\sqrt{\pi}}\left(
\frac13-\frac{c^{2}}6
\right)
+\erfc\left(\frac{c}{2}\right)\,\frac{c^{3}}{12}
\Bigg]\\
A^{(2)}_{4}
&=
\sigma_{1}\Bigg[
\frac{e^{-c^{2}/4}}{\sqrt{\pi}}\Bigg(
-\frac{\alpha_{2}\beta_{2}^{2}}{24}c
+\frac{\alpha_{2}\beta_{2}^{2}}{24}c^{3}
-\frac{\beta_{3}k}{6}c
+\frac{\beta_{3}k}{12}c^{3}
-\frac{5\beta_{2}k^{2}}{24}c
+\frac{7\beta_{2}k^{2}}{24}c^{3}
\Bigg)\notag\\
&\qquad
+\erfc\left(\frac{c}{2}\right)\Bigg(
-\frac{\alpha_{2}\beta_{2}^{2}}{48}c^{4}
-\frac{\beta_{3}k}{24}c^{4}
-\frac{7\beta_{2}k^{2}}{48}c^{4}
\Bigg)
\Bigg]\notag\\
&\quad+
\sigma_{2}\Bigg[
\frac{e^{-c^{2}/4}}{\sqrt{\pi}}\Bigg(
-\frac{\alpha_{3}\beta_{2}}{12}c
-\frac{\alpha_{2}\beta_{3}}{12}c
+\frac{\alpha_{3}\beta_{2}}{24}c^{3}
+\frac{\alpha_{2}\beta_{3}}{24}c^{3}
+\frac{\alpha_{1}\beta_{1}k}{12}c
-\frac{\alpha_{2}\beta_{2}k}{8}c
-\frac{\alpha_{1}\beta_{1}k}{24}c^{3}\notag\\
&\qquad
+\frac{\alpha_{2}\beta_{2}k}{8}c^{3}
-\frac{k^{3}}{8}c
+\frac{3k^{3}}{8}c^{3}
\Bigg)
+\erfc\left(\frac{c}{2}\right)\Bigg(
-\frac{\alpha_{3}\beta_{2}}{48}c^{4}
-\frac{\alpha_{2}\beta_{3}}{48}c^{4}
+\frac{\alpha_{1}\beta_{1}k}{48}c^{4}
-\frac{\alpha_{2}\beta_{2}k}{16}c^{4}
-\frac{3k^{3}}{16}c^{4}
\Bigg)
\Bigg]\notag\\
&\quad+
\sigma_{1}'\Bigg[
\frac{e^{-c^{2}/4}}{\sqrt{\pi}}\left(
-\frac{\beta_{1}k}{2}c
+\frac{\beta_{1}k}{4}c^{3}
\right)
-\erfc\left(\frac{c}{2}\right)\frac{\beta_{1}k}{8}c^{4}
\Bigg]\notag\\
&\quad+
\sigma_{2}'\Bigg[
\frac{e^{-c^{2}/4}}{\sqrt{\pi}}\left(
-\frac{\beta_{1}\alpha_{2}}{12}c
-\frac{\alpha_{1}\beta_{2}}{12}c
+\frac{\beta_{1}\alpha_{2}}{24}c^{3}
+\frac{\alpha_{1}\beta_{2}}{24}c^{3}
\right)
-\erfc\left(\frac{c}{2}\right)\left(
\frac{\beta_{1}\alpha_{2}}{48}c^{4}
+\frac{\alpha_{1}\beta_{2}}{48}c^{4}
\right)
\Bigg]\notag\\
&\quad+
\sigma_{1}''\Bigg[
\frac{e^{-c^{2}/4}}{\sqrt{\pi}}\left(
\frac{\beta_{2}}{6}c-\frac{\beta_{2}}{12}c^{3}
\right)
+\erfc\left(\frac{c}{2}\right)\,\frac{\beta_{2}}{24}c^{4}
\Bigg]\notag\\
&\quad+
\sigma_{2}''\Bigg[
\frac{e^{-c^{2}/4}}{\sqrt{\pi}}\left(
\frac{k}{2}c-\frac{k}{4}c^{3}
\right)
+\erfc\left(\frac{c}{2}\right)\,\frac{k}{8}c^{4}
\Bigg]\notag\\
&\quad+
\sigma_{1}\,\gamma^{(3)}\Bigg[
\frac{e^{-c^{2}/4}}{\sqrt{\pi}}\left(
-\frac{\beta_{1}}{6}c+\frac{\beta_{1}}{12}c^{3}
\right)
-\erfc\left(\frac{c}{2}\right)\,\frac{\beta_{1}}{24}c^{4}
\Bigg]\notag\\
&\quad+
\sigma_{2}'\,\gamma^{(3)}\Bigg[
\frac{e^{-c^{2}/4}}{\sqrt{\pi}}\left(
\frac13 c-\frac16 c^{3}
\right)
+\erfc\left(\frac{c}{2}\right)\,\frac1{12}c^{4}
\Bigg]\notag\\
&\quad+
\sigma_{2}\,\gamma^{(4)}\Bigg[
\frac{e^{-c^{2}/4}}{\sqrt{\pi}}\left(
\frac1{12}c-\frac1{24}c^{3}
\right)
+\erfc\left(\frac{c}{2}\right)\,\frac1{48}c^{4}
\Bigg].
\end{align}
The expansion shows that the leading-order behavior is identical to the scalar case, with coupling entering only at the next order through the coefficients $\alpha_2$ and $\beta_2$. Higher-order terms reflect the full structure of the system symbol, with mixed coefficients $\alpha_j \beta_j$, $j=1,2,3$, appearing at order $\veps^{3/2}$ and beyond. As for the volume potential, coupling contributions enter at the stage where interactions between components are resolved by the symbol, and thereafter propagate systematically through the expansion.
\subsection{The Double Layer Potential $\dlp$}\label{subsec:dlpsys}
The treatment of the vector-valued double layer potential $\bm{\calD}[\vmu]$ 
is analogous to \cref{subsec:slpsys}, with density $\vmu = [\mu_1\;\mu_2]$ 
and source $\vF = -\nabla\cdot(\vnu\vmu\delta(g))$ interpreted componentwise. The additional factor of 
$i\vnu\cdot\veta$ modifies only the algebraic structure of the intermediate 
expansions. The asymptotic expansion takes the form
\begin{equation}
\bm{\calD}_L[\vmu](\vx) = \sum_{p=0}^{P}\bm{A}_p\veps^{p/2} + 
\ordo{\veps}{(P+1)/2},
\end{equation}
with $\bm{A}_p = [A^{(1)}_p\;A^{(2)}_p]$. For $P = 3$ we have
\begin{align}
A^{(1)}_{0}
&=
-\frac12\,\erfc\left(\frac{c}{2}\right)\mu_{1}\\
A^{(1)}_{1}
&=
-\mu_{1}\,\frac{e^{-c^{2}/4}}{\sqrt{\pi}}\,\frac{k}{2}
+\mu_{2}\left(
\frac{\alpha_{2}e^{-c^{2}/4}}{\sqrt{\pi}}\,\frac12
-\frac{\alpha_{2}}{2}\,c\,\erfc\left(\frac{c}{2}\right)
\right)\\
A^{(1)}_{2}
&=
-\mu_{1}\Bigg[
\frac{e^{-c^{2}/4}}{\sqrt{\pi}}\left(
-\frac{3\alpha_{2}\beta_{2}}{8}c
+\frac{3k^{2}}{8}c
\right)
+\erfc\left(\frac{c}{2}\right)\frac{\alpha_{2}\beta_{2}}{4}c^{2}
\Bigg]\notag\\
&\quad
-\mu_{2}\Bigg[
\frac{e^{-c^{2}/4}}{\sqrt{\pi}}\left(
-\frac{\alpha_{3}}{2}c-\frac{\alpha_{2}k}{4}c
\right)
+\erfc\left(\frac{c}{2}\right)\left(
\frac{\alpha_{3}}{4}c^{2}+\frac{\alpha_{2}k}{4}c^{2}
\right)
\Bigg]\notag\\
&\quad
-\mu_{2}'\Bigg[
-\frac{\alpha_{1}}{2}\frac{e^{-c^{2}/4}}{\sqrt{\pi}}\,c
+\frac{\alpha_{1}}{4}\,c^{2}\erfc\left(\frac{c}{2}\right)
\Bigg]
-\mu_{1}''\Bigg[
\frac{e^{-c^{2}/4}}{\sqrt{\pi}}\frac{c}{2}
-\erfc\left(\frac{c}{2}\right)\frac{c^{2}}{4}
\Bigg]\\
A^{(1)}_{3}
&=
-\mu_{1}\Bigg[
\frac{e^{-c^{2}/4}}{\sqrt{\pi}}\Bigg(
\frac{\alpha_{3}\beta_{2}}{12}
+\frac{\alpha_{2}\beta_{3}}{12}
-\frac{\alpha_{3}\beta_{2}}{6}c^{2}
-\frac{\alpha_{2}\beta_{3}}{6}c^{2}
-\frac{\alpha_{1}\beta_{1}k}{12}
+\frac{\alpha_{1}\beta_{1}k}{6}c^{2}
+\frac{\alpha_{2}\beta_{2}k}{24}
-\frac{13\alpha_{2}\beta_{2}k}{48}c^{2}\notag\\
&\qquad
-\frac{5k^{3}}{8}
+\frac{5k^{3}}{16}c^{2}
\Bigg)
+\erfc\left(\frac{c}{2}\right)\left(
\frac{\alpha_{3}\beta_{2}}{12}c^{3}
+\frac{\alpha_{2}\beta_{3}}{12}c^{3}
-\frac{\alpha_{1}\beta_{1}k}{12}c^{3}
+\frac{\alpha_{2}\beta_{2}k}{6}c^{3}
\right)
\Bigg]\notag\\
&\quad
-\mu_{2}\Bigg[
\frac{e^{-c^{2}/4}}{\sqrt{\pi}}\left(
\frac{\alpha_{2}^{2}\beta_{2}}{24}
-\frac{7\alpha_{2}^{2}\beta_{2}}{48}c^{2}
-\frac{\alpha_{3}k}{6}
-\frac{\alpha_{3}k}{6}c^{2}
-\frac{\alpha_{2}k^{2}}{8}
-\frac{5\alpha_{2}k^{2}}{16}c^{2}
\right)\notag\\
&\qquad
+\erfc\left(\frac{c}{2}\right)\left(
\frac{\alpha_{2}^{2}\beta_{2}}{12}c^{3}
+\frac{\alpha_{3}k}{12}c^{3}
+\frac{\alpha_{2}k^{2}}{4}c^{3}
\right)
\Bigg]\notag\\
&\quad
-\mu_{1}'\Bigg[
\frac{e^{-c^{2}/4}}{\sqrt{\pi}}\left(
\frac{\alpha_{2}\beta_{1}}{12}
+\frac{\alpha_{1}\beta_{2}}{12}
-\frac{\alpha_{2}\beta_{1}}{6}c^{2}
-\frac{\alpha_{1}\beta_{2}}{6}c^{2}
\right)
+\erfc\left(\frac{c}{2}\right)\left(
\frac{\alpha_{2}\beta_{1}}{12}c^{3}
+\frac{\alpha_{1}\beta_{2}}{12}c^{3}
\right)
\Bigg]\notag\\
&\quad
-\mu_{2}'\Bigg[
\frac{\alpha_{1}k\,e^{-c^{2}/4}}{\sqrt{\pi}}\left(
-\frac16-\frac23 c^{2}
\right)
+\erfc\left(\frac{c}{2}\right)\frac{\alpha_{1}k}{3}c^{3}
\Bigg]\notag\\
&\quad
-\mu_{1}''\Bigg[
\frac{k\,e^{-c^{2}/4}}{\sqrt{\pi}}\left(
\frac12+\frac12 c^{2}
\right)
-\erfc\left(\frac{c}{2}\right)\frac{k}{4}c^{3}
\Bigg]\notag\\
&\quad
-\mu_{2}''\Bigg[
\frac{\alpha_{2}e^{-c^{2}/4}}{\sqrt{\pi}}\left(
-\frac16+\frac13 c^{2}
\right)
-\erfc\left(\frac{c}{2}\right)\frac{\alpha_{2}}{6}c^{3}
\Bigg]\notag\\
&\quad
-\mu_{2}\gamma^{(3)}\Bigg[
\frac{\alpha_{1}e^{-c^{2}/4}}{\sqrt{\pi}}\left(
-\frac16-\frac16 c^{2}
\right)
+\erfc\left(\frac{c}{2}\right)\frac{\alpha_{1}}{12}c^{3}
\Bigg]\notag\\
&\quad
-\mu_{1}'\gamma^{(3)}\Bigg[
\frac{e^{-c^{2}/4}}{\sqrt{\pi}}\left(
\frac23+\frac16 c^{2}
\right)
-\erfc\left(\frac{c}{2}\right)\frac{1}{12}c^{3}
\Bigg]
-\mu_{1}\gamma^{(4)}\Bigg[
\frac{e^{-c^{2}/4}}{\sqrt{\pi}}\frac14
\Bigg]
\end{align}
and
\begin{align}
A^{(2)}_{0}
&=
-\frac12\,\erfc\left(\frac{c}{2}\right)\mu_{2}\\
A^{(2)}_{1}
&=
-\mu_{1}\left(
-\frac{\beta_{2}e^{-c^{2}/4}}{2\sqrt{\pi}}
+\frac{\beta_{2}c}{2}\,\erfc\left(\frac{c}{2}\right)
\right)
-\mu_{2}\,\frac{k\,e^{-c^{2}/4}}{2\sqrt{\pi}}\\
A^{(2)}_{2}
&=
-\mu_{1}\Bigg[
\frac{c\,e^{-c^{2}/4}}{\sqrt{\pi}}\left(
-\frac{\beta_{3}}{2}
-\frac{\beta_{2}k}{4}
\right)
+\left(
\frac{\beta_{3}c^{2}}{4}
+\frac{\beta_{2}kc^{2}}{4}
\right)\erfc\left(\frac{c}{2}\right)
\Bigg]\notag\\
&\quad
-\mu_{2}\Bigg[
\frac{c\,e^{-c^{2}/4}}{\sqrt{\pi}}\left(
-\frac{3\alpha_{2}\beta_{2}}{8}
+\frac{3k^{2}}{8}
\right)
+\frac{\alpha_{2}\beta_{2}c^{2}}{4}\,\erfc\left(\frac{c}{2}\right)
\Bigg]\notag\\
&\quad
-\mu_{1}'\left[
-\frac{\beta_{1}c\,e^{-c^{2}/4}}{2\sqrt{\pi}}
+\frac{\beta_{1}c^{2}}{4}\,\erfc\left(\frac{c}{2}\right)
\right]
-\mu_{2}''\left[
\frac{c\,e^{-c^{2}/4}}{2\sqrt{\pi}}
-\frac{c^{2}}{4}\,\erfc\left(\frac{c}{2}\right)
\right]\\
A^{(2)}_{3}
&=
-\mu_{1}\Bigg[
\frac{e^{-c^{2}/4}}{\sqrt{\pi}}\left(
\frac{\alpha_{2}\beta_{2}^{2}}{24}
-\frac{7\alpha_{2}\beta_{2}^{2}}{48}c^{2}
-\frac{\beta_{3}k}{6}
-\frac{\beta_{3}k}{6}c^{2}
-\frac{\beta_{2}k^{2}}{8}
-\frac{5\beta_{2}k^{2}}{16}c^{2}
\right)\notag\\
&\qquad
+\left(
\frac{\alpha_{2}\beta_{2}^{2}}{12}c^{3}
+\frac{\beta_{3}k}{12}c^{3}
+\frac{\beta_{2}k^{2}}{4}c^{3}
\right)\erfc\left(\frac{c}{2}\right)
\Bigg]\notag\\
&\quad
-\mu_{2}\Bigg[
\frac{e^{-c^{2}/4}}{\sqrt{\pi}}\Bigg(
\frac{\alpha_{3}\beta_{2}+\alpha_{2}\beta_{3}}{12}
-\frac{\alpha_{3}\beta_{2}+\alpha_{2}\beta_{3}}{6}c^{2}
-\frac{\alpha_{1}\beta_{1}k}{12}
+\frac{\alpha_{1}\beta_{1}k}{6}c^{2}
+\frac{\alpha_{2}\beta_{2}k}{24}
-\frac{13\alpha_{2}\beta_{2}k}{48}c^{2}\notag\\
&\qquad
-\frac{5k^{3}}{8}
+\frac{5k^{3}}{16}c^{2}
\Bigg)
+\left(
\frac{\alpha_{3}\beta_{2}+\alpha_{2}\beta_{3}}{12}c^{3}
-\frac{\alpha_{1}\beta_{1}k}{12}c^{3}
+\frac{\alpha_{2}\beta_{2}k}{6}c^{3}
\right)\erfc\left(\frac{c}{2}\right)
\Bigg]\notag\\
&\quad
-\mu_{1}'\Bigg[
\frac{\beta_{1}k\,e^{-c^{2}/4}}{\sqrt{\pi}}\left(
-\frac16-\frac23 c^{2}
\right)
+\frac{\beta_{1}k}{3}c^{3}\,\erfc\left(\frac{c}{2}\right)
\Bigg]\notag\\
&\quad
-\mu_{2}'\Bigg[
\frac{e^{-c^{2}/4}}{\sqrt{\pi}}\left(
\frac{\alpha_{2}\beta_{1}+\alpha_{1}\beta_{2}}{12}
-\frac{\alpha_{2}\beta_{1}+\alpha_{1}\beta_{2}}{6}c^{2}
\right)
+\frac{\alpha_{2}\beta_{1}+\alpha_{1}\beta_{2}}{12}c^{3}\,\erfc\left(\frac{c}{2}\right)
\Bigg]\notag\\
&\quad
-\mu_{1}''\Bigg[
\frac{\beta_{2}e^{-c^{2}/4}}{\sqrt{\pi}}\left(
-\frac16+\frac13 c^{2}
\right)
-\frac{\beta_{2}}{6}c^{3}\,\erfc\left(\frac{c}{2}\right)
\Bigg]\notag\\
&\quad
-\mu_{2}''\Bigg[
\frac{k\,e^{-c^{2}/4}}{\sqrt{\pi}}\left(
\frac12+\frac12 c^{2}
\right)
-\frac{k}{4}c^{3}\,\erfc\left(\frac{c}{2}\right)
\Bigg]\notag\\
&\quad
-\mu_{1}\gamma^{(3)}\Bigg[
\frac{\beta_{1}e^{-c^{2}/4}}{\sqrt{\pi}}\left(
-\frac16-\frac16 c^{2}
\right)
+\frac{\beta_{1}}{12}c^{3}\,\erfc\left(\frac{c}{2}\right)
\Bigg]\notag\\
&\quad
-\mu_{2}'\gamma^{(3)}\Bigg[
\frac{e^{-c^{2}/4}}{\sqrt{\pi}}\left(
\frac23+\frac16 c^{2}
\right)
-\frac{c^{3}}{12}\,\erfc\left(\frac{c}{2}\right)
\Bigg]
-\mu_{2}\gamma^{(4)}\,
\frac{e^{-c^{2}/4}}{4\sqrt{\pi}}.
\end{align}
We see that coupling appears already at order $\sqrt{\veps}$, reflecting the presence of first-order terms in the operator and their interaction with the normal direction. At higher orders, the expansion incorporates both the system coupling and geometric corrections in a fully coupled manner. As for the volume potential, coupling terms enter at the stage where derivatives interact with the symbol, and subsequently influence all higher-order contributions.
\section{Numerical results}\label{sec:numerical_results}
In these numerical experiments, we verify the asymptotic formulas derived for the Poisson equation in $\mathbb{R}^3$ in \cref{sec:poissonr3} and for coupled systems in \cref{sec:coupsys}. We do not solve the PDEs directly, but instead compute the solutions using Green's third identity. We therefore assume that all required data are available analytically, including high-order derivatives of the solution and of the surface.

Given quadrature points and weights ${(\vx_i,w_i)}_{i=1}^N$ for the volume mesh, we follow \cite{lghtwght_ps2026} and set $\delta = 3h^2$, where $h = \left(\frac{1}{N}\sum_{i=1}^N w_i\right)^{1/d}$ is the average mesh resolution and $d$ is the dimension. The history parts are computed using the non-uniform FFT library \cite{finufftlib}.

\subsection{The Poisson equation in a torus}
We validate the asymptotic formulas derived for the three-dimensional Poisson equation by evaluating the solution using Green’s representation formula in a nontrivial geometry and measuring the resulting convergence behavior. To this end, we consider the manufactured solution
\begin{equation}
  \begin{split}\label{eq:linsys_sol}
u(x,y,z)
&= \exp\big((x+2) + 2(y+2) + 3(z+1)\big)
+ \sin\big((x+2)^2 (y+2) - (y+2)(z+1)^2 + (z+1)\big)\\
&+ (x+2)(y+2)^2 (z+1)^3
    - \cos\big((x+2)(y+2)(z+1)\big),
      \end{split}
\end{equation}
from which the corresponding volume and boundary data are obtained analytically.

We evaluate the solution $u$ in the interior of a torus with major radius $0.5$ and minor radius $0.3$ at $1000$ target points sampled randomly in a thin interior layer adjacent to the boundary, within a distance of at most $0.03$ from $\partial\Omega$. The solution is then computed using Green’s representation formula, introduced in \cref{sec:pottheory}, which we restate here in operator form:
\begin{equation}
u(\vx)
=
-\vp[\Delta u](\vx)
+
\slp\!\left[\frac{\partial u}{\partial \vnu}\right](\vx)
-
\dlp[u](\vx).
\end{equation}
Here, the single-layer density is given by the normal derivative $\partial u/\partial \vnu$, while the double-layer density is given by the trace of $u$ on $\partial\Omega$.

The volume, single layer, and double layer contributions are evaluated using the local--history decomposition described in the preceding sections. The volumetric discretization is based on a tensor-product quadrature, with Gauss--Legendre nodes in the radial direction and periodic trapezoidal rules in the angular directions, yielding high-order accuracy for smooth integrands. The surface discretization is obtained using the library \texttt{fmm3dbie}~\cite{fmm3dbie}, based on patchwise parameterizations with fourth-order Vioreanu--Rokhlin nodes.

\begin{figure}[t]
\centering
  \begin{tikzpicture}[scale=0.78]
\definecolor{mycolor1}{rgb}{0.00000,0.44700,0.74100}%
\definecolor{mycolor2}{rgb}{0.85000,0.32500,0.09800}%
\begin{axis}[%
scale only axis,
xmode=log,
xmin=0.007,
xmax=0.05,
xminorticks=true,
xlabel style={font=\Large\color{white!15!black}},
xlabel={$h$},
xticklabel style={font=\large},
ymode=log,
ymin=8e-5,
ymax=3e-1,
yminorticks=true,
ylabel style={font=\Large\color{white!15!black}},
ylabel={Relative max error},
yticklabel style={font=\large},
axis background/.style={fill=white},
xmajorgrids,
xminorgrids,
ymajorgrids,
yminorgrids,
legend style={font=\large,at={(0.03,0.97)}, anchor=north west, legend cell align=left, align=left, draw=white!15!black},
clip mode=individual,
]

\addplot [color=mycolor1, line width=1.0pt, mark=*, mark options={solid, mycolor1}, mark size=2.0pt]
  table[row sep=crcr]{%
     4.936536597953740e-02     1.888942787154664e-01 \\
     4.442882938158366e-02     1.493622181226736e-01 \\
     3.417602260121820e-02     7.978125032027909e-02 \\
     3.173487812970262e-02     6.614995553004702e-02 \\
     2.776801836348979e-02     4.672018660425230e-02 \\
     2.338359441135982e-02     2.935035653254301e-02 \\
     2.221441469079183e-02     2.546015318764789e-02 \\
     1.851201224232652e-02     1.519078171250108e-02 \\
     1.708801130060910e-02     1.213752951199591e-02 \\
     1.480960979386122e-02     8.087514609633941e-03 \\
     1.269395125188105e-02     5.191257739730585e-03 \\
     1.083629984916675e-02     3.274946860912651e-03 \\
     9.658441169909491e-03     2.335086995940257e-03 \\
     8.382797996525219e-03     1.534174662125217e-03 \\
};
\addlegendentry{error ($P=2$)}

\addplot [color=mycolor2, dotted, line width=1.0pt, mark=o, mark options={solid, mycolor2}, mark size=2.0pt]
  table[row sep=crcr]{%
     4.936536597953740e-02     1.633416326247415e-01 \\
     4.442882938158366e-02     1.080922816651702e-01 \\
     3.417602260121820e-02     3.772023296230041e-02 \\
     3.173487812970262e-02     2.788192827547676e-02 \\
     2.776801836348979e-02     1.611698633890598e-02 \\
     2.338359441135982e-02     7.929098538088895e-03 \\
     2.221441469079183e-02     6.414116715716092e-03 \\
     1.851201224232652e-02     3.020860831829677e-03 \\
     1.708801130060910e-02     2.185696546150119e-03 \\
     1.480960979386122e-02     1.228866706843361e-03 \\
     1.269395125188105e-02     6.635760005370665e-04 \\
     1.083629984916675e-02     3.542957616701091e-04 \\
     9.658441169909491e-03     2.252800124398748e-04 \\
     8.382797996525219e-03     1.323911977160450e-04 \\
};
\addlegendentry{error ($P=3$)}

\addplot [color=black, dashed, line width=1.0pt]
  table[row sep=crcr]{%
     2.338359441135982e-02     1.974450613217327e-02 \\
     2.221441469079183e-02     1.692844594507206e-02 \\
     1.851201224232652e-02     9.796554366361144e-03 \\
     1.708801130060910e-02     7.705255323200755e-03 \\
     1.480960979386122e-02     5.015835835576907e-03 \\
     1.269395125188105e-02     3.158660467885166e-03 \\
     1.083629984916675e-02     1.964968116547590e-03 \\
     9.658441169909491e-03     1.391340999841543e-03 \\
     8.382797996525219e-03     9.096607774241589e-04 \\
};
\addlegendentry{$O(h^{3})$}

\addplot [color=black, dashdotted, line width=1.0pt]
  table[row sep=crcr]{%
     3.417602260121820e-02     2.103193471781847e-02 \\
     3.173487812970262e-02     1.563653393053972e-02 \\
     2.776801836348979e-02     9.165849113092248e-03 \\
     2.338359441135982e-02     4.609334546816043e-03 \\
     2.221441469079183e-02     3.754331796722586e-03 \\
     1.851201224232652e-02     1.810538096413283e-03 \\
     1.708801130060910e-02     1.314495919863655e-03 \\
     1.480960979386122e-02     7.415964042908809e-04 \\
     1.269395125188105e-02     4.002952686218166e-04 \\
     1.083629984916675e-02     2.125774570020655e-04 \\
};
\addlegendentry{$O(h^{4})$}

\end{axis}
\end{tikzpicture}
\caption{Convergence of the numerical method. For truncation levels $P=2$ and $P=3$ in the asymptotic expansions, the errors display third- and fourth-order convergence, respectively, consistent with the predicted $\mathcal{O}(h^{(P+1)})$ scaling. The dashed lines denote the reference slopes $\ordo{h}{3}$ and $\ordo{h}{4}$.}
\label{fig:convergence}
\end{figure}

The volume discretization is coupled to the surface resolution by selecting the number of volume nodes so that the corresponding volumetric spacing satisfies $h_v \approx h$, where $h \sim \sqrt{\mathrm{area}/N_s}$ denotes the average spacing between surface nodes. This ensures that both discretizations resolve the same physical scales and avoids imbalances in the resulting approximation errors.

Errors are measured in the relative maximum norm over all target points, and we report results for both $P=2$ and $P=3$ in the asymptotic formulas. The convergence results are shown in \cref{fig:convergence}, where we observe the expected orders of accuracy predicted by the asymptotic expansions derived in \cref{sec:poissonr3}. In particular, the results confirm that the local asymptotic corrections recover the singular behavior near the boundary, leading to high-order convergence under refinement.

\subsection{A coupled strongly elliptic system}
\begin{figure}[t]
\centering
  \includegraphics[width=0.4\linewidth]{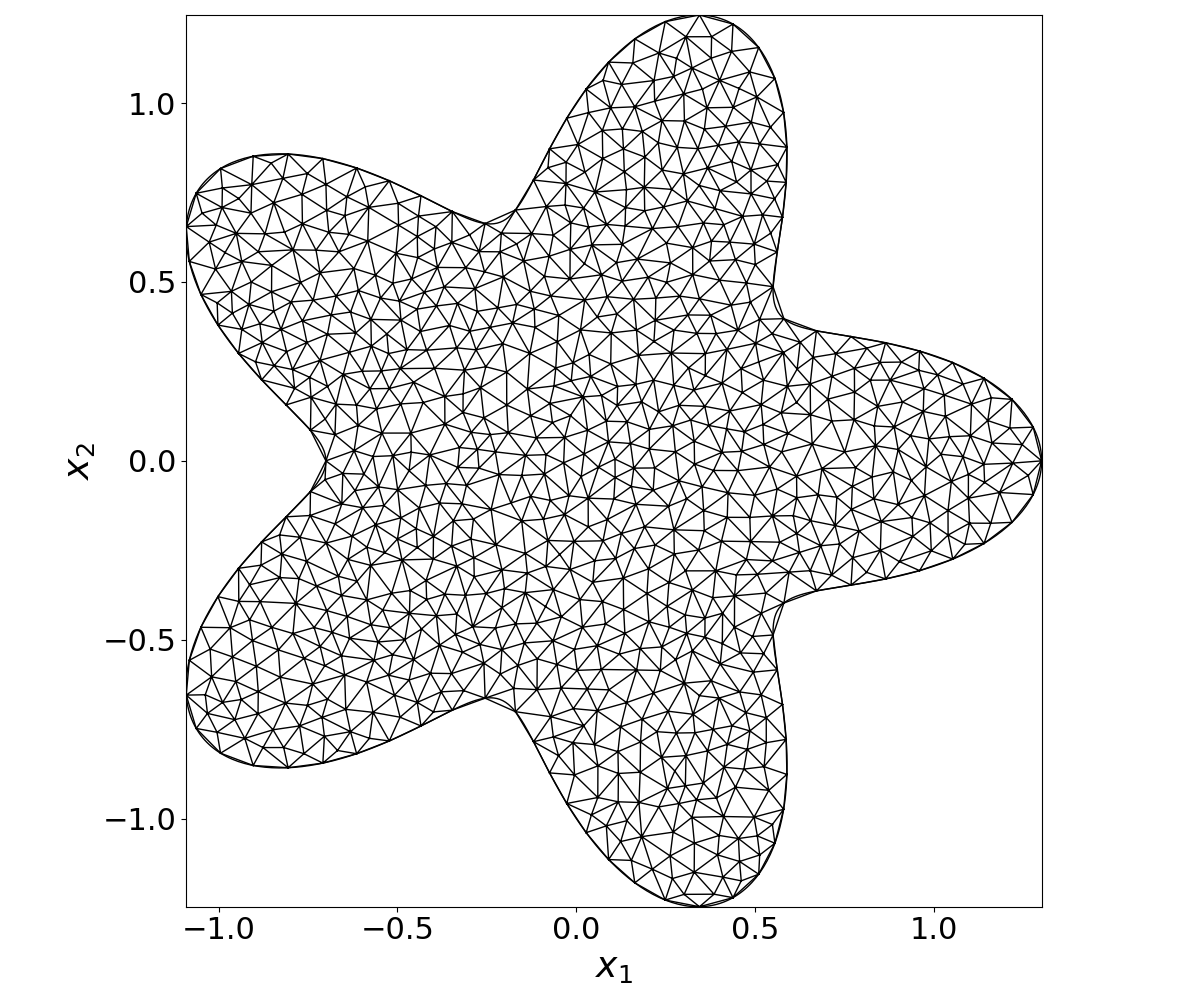}
\caption{A volumetric mesh for the domain $\Omega$ defined by the parametrization \eqref{eq:bdryparam} of the boundary curve.}
\label{fig:starmesh}
\end{figure}

We verify the asymptotic formulas for $\slp$, $\dlp$ and $\vp$ for coupled linear system of PDEs \eqref{eq:coupsys} by computing the solution $\vv$ using Green's third identity
  \begin{equation}\label{eq:sysgreen3rd}
          \vv
          =
          \bm{\vp}[\vf]
          +
          \bm{\slp}[\partial_{\vnu} \vv - D_{\vnu} \vv]
          -
          \bm{\dlp}[\vv]
        \end{equation}
        where
        \begin{equation}
          \partial_{\vnu} \vv - D_{\vnu} \vv
          =
          \begin{bmatrix}
\partial_{\vnu} v_1 - (\va\cdot\vnu)v_2\\
\partial_{\vnu} v_2 - (\vb\cdot\vnu)v_1
          \end{bmatrix}.
        \end{equation}
        The solution is given by
        \begin{equation}
\vv(\vx)
=\begin{bmatrix}
\sin(2.3\pi x_1+0.4)\cos(1.7\pi x_2-0.2)+0.7x_1x_2-0.4(x_1-0.3)^2(x_2+0.6)\\
\cos(1.6\pi x_1-0.1)\sin(2.1\pi x_2+0.3)-0.5x_1^2x_2+0.8(x_1+0.2)(x_2-0.4)^2
\end{bmatrix},
\end{equation}
and the boundary curve is described by
\begin{equation}\label{eq:bdryparam}
  \vga(s) = \left(1+0.3\cos(5s)\right) \begin{bmatrix}\cos(s)\\\sin(s)\end{bmatrix},\quad {\rm for}\ s\in[0,2\pi),
\end{equation}
and shown in \cref{fig:starmesh}. The boundary is discretized into chunks of equal length in parameter $s$. We generate a $p$th order triangle mesh using $p$th order boundary elements and a $p$th order quadrature rule on each triangle. The meshes are generated using {\tt MeshPy} \cite{meshpylib} as an interface to {\tt Triangle} \cite{trianglepaper,trianglelib}. As quadrature rule over each triangle element we use Vioreanu-Rokhlin nodes \cite{vioreanu2014sisc}. We generate a sequence of 6 meshes, corresponding to $h$ approximately equal to $0.007$, $0.006$, $0.004$, $0.002$, $0.001$, and $0.0007$. The representation \cref{eq:sysgreen3rd} is evaluated at the quadrature nodes for each mesh, and the errors are measured in the relative discrete $\ell^{\infty}$ and $\ell^2$ norms. In the asymptotic expansions we keep terms up to and including $\delta^2$, which translates to a $\ordo{h}{5}$ method. To have an overall fifth order method, we set $p=5$.

First, we set $\va = [2,-1.31]$, $\vb = [-7.23,4.3]$, $a_3=1$ and $b_3=-1$. This provides a baseline test in which the first-order coupling terms are of moderate size. The left panel of \cref{fig:5thorder} shows the expected $\ordo{h}{5}$ convergence in both the relative maximum norm and the relative discrete $\ell^2$ norm, confirming that the terms retained through $\delta^2$ give the predicted fifth-order accuracy when combined with fifth-order volume quadrature.

We then repeat the experiment with $a_1=200$, keeping all other parameters unchanged. This is a more demanding test, since the first-order coupling terms are now large and the different contributions in the representation formula have significantly different magnitudes. The right panel of \cref{fig:5thorder} again shows fifth-order convergence. The main visible effect is that the error for $v_1$ is shifted upward by roughly one digit, which reflects increased numerical sensitivity rather than a change in the observed convergence order.

\begin{figure}[ht]
\centering

\begin{subfigure}{0.49\textwidth}
  \centering
\begin{tikzpicture}[scale=0.62]

\begin{axis}[
scale only axis,
xmode=log,
ymode=log,
xmin=6e-4,
xmax=8e-3,
xminorticks=true,
xlabel={$h$},
xticklabel style={font=\large},
ymin=1e-10,
ymax=1e-3,
yminorticks=true,
ylabel={Relative error},
yticklabel style={font=\large},
xlabel style={font=\Large\color{white!15!black}},
ylabel style={font=\Large\color{white!15!black}},
axis background/.style={fill=white},
xmajorgrids,
xminorgrids,
ymajorgrids,
yminorgrids,
legend style={
    font=\Large,
    at={(0.63,0.37)},
    anchor=north west,
    legend cell align=left,
    align=left,
    draw=white!15!black
},
clip mode=individual,
]

\addplot [color=blue, line width=1.0pt, mark=*, mark options={solid}, mark size=2.0pt]
table[row sep=crcr]{
7.230824819348436e-03  3.303990623141813e-04 \\
5.787060338230405e-03  1.187848413755858e-04 \\
3.877843685051529e-03  1.782539188380658e-05 \\
2.103067914659541e-03  9.008909653075288e-07 \\
9.087673506033234e-04  1.397956095773128e-08 \\
6.907543842015960e-04  3.551587382606267e-09 \\
};
\addlegendentry{$u$, $e_{\infty}$}

\addplot [color=red,  line width=1.0pt, mark=o, mark options={solid}, mark size=2.0pt]
table[row sep=crcr]{
7.230824819348436e-03  5.120258231597507e-05 \\
5.787060338230405e-03  1.724194778988068e-05 \\
3.877843685051529e-03  2.571969474299824e-06 \\
2.103067914659541e-03  1.176254183045422e-07 \\
9.087673506033234e-04  1.444954959110085e-09 \\
6.907543842015960e-04  4.177374123751055e-10 \\
};
\addlegendentry{$u$, $e_{2}$}

\addplot [color=magenta, dashed,line width=1.0pt, mark=square*, mark options={solid}, mark size=2.0pt]
table[row sep=crcr]{
7.230824819348436e-03  2.425714893753772e-04 \\
5.787060338230405e-03  8.759514341383189e-05 \\
3.877843685051529e-03  1.324529210488022e-05 \\
2.103067914659541e-03  6.683928971328392e-07 \\
9.087673506033234e-04  1.035436758921886e-08 \\
6.907543842015960e-04  2.634274447639177e-09 \\
};
\addlegendentry{$v$, $e_{\infty}$}

\addplot [color=teal, dashed, line width=1.0pt, mark=diamond*, mark options={solid}, mark size=2.0pt]
table[row sep=crcr]{
7.230824819348436e-03  4.653095499224069e-05 \\
5.787060338230405e-03  1.518438706730287e-05 \\
3.877843685051529e-03  2.329544957612467e-06 \\
2.103067914659541e-03  1.055394666360268e-07 \\
9.087673506033234e-04  1.330627621039948e-09 \\
6.907543842015960e-04  4.062096684598654e-10 \\
};
\addlegendentry{$v$, $e_{2}$}

\addplot [color=black, dotted, line width=1.0pt]
table[row sep=crcr]{
7.230824819348436e-03  1.300728899471055e-04 \\
5.787060338230405e-03  4.270694506668837e-05 \\
3.877843685051529e-03  5.768338791095939e-06 \\
2.103067914659541e-03  2.706955062428745e-07 \\
9.087673506033234e-04  4.077342769325366e-09 \\
6.907543842015960e-04  1.034667101051845e-09 \\
};
\addlegendentry{$O(h^{5})$}

\end{axis}
\end{tikzpicture}

\end{subfigure}
\hfill
\begin{subfigure}{0.49\textwidth}
  \centering
  \begin{tikzpicture}[scale=0.62]

\begin{axis}[
scale only axis,
xmode=log,
ymode=log,
xmin=6e-4,
xmax=8e-3,
xminorticks=true,
xlabel={$h$},
xticklabel style={font=\large},
ymin=1e-10,
ymax=1e-3,
yminorticks=true,
ylabel={Relative error},
yticklabel style={font=\large},
xlabel style={font=\Large\color{white!15!black}},
ylabel style={font=\Large\color{white!15!black}},
axis background/.style={fill=white},
xmajorgrids,
xminorgrids,
ymajorgrids,
yminorgrids,
legend style={
    font=\Large,
    at={(0.63,0.37)},
    anchor=north west,
    legend cell align=left,
    align=left,
    draw=white!15!black
},
clip mode=individual,
]

\addplot [color=blue, line width=1.0pt, mark=*, mark options={solid}, mark size=2.0pt]
table[row sep=crcr]{
7.230824819348436e-03  2.386608008281675e-03 \\
5.787060338230405e-03  8.543600610971848e-04 \\
3.877843685051529e-03  1.273547427162854e-04 \\
2.103067914659541e-03  6.377192154753944e-06 \\
9.087673506033234e-04  9.836259695927087e-08 \\
6.907543842015960e-04  2.987070173634372e-08 \\
};
\addlegendentry{$u$, $e_{\infty}$}

\addplot [color=red,  line width=1.0pt, mark=o, mark options={solid}, mark size=2.0pt]
table[row sep=crcr]{
7.230824819348436e-03  4.076736655780911e-04 \\
5.787060338230405e-03  1.335837454154048e-04 \\
3.877843685051529e-03  1.981810013676633e-05 \\
2.103067914659541e-03  8.793192620108784e-07 \\
9.087673506033234e-04  1.240161244810811e-08 \\
6.907543842015960e-04  4.520133772959494e-09 \\
};
\addlegendentry{$u$, $e_{2}$}

\addplot [color=magenta, dashed,line width=1.0pt, mark=square*, mark options={solid}, mark size=2.0pt]
table[row sep=crcr]{
7.230824819348436e-03  3.591681085367163e-04 \\
5.787060338230405e-03  1.286030565311815e-04 \\
3.877843685051529e-03  1.925501808079003e-05 \\
2.103067914659541e-03  9.641587126687263e-07 \\
9.087673506033234e-04  1.491036184109868e-08 \\
6.907543842015960e-04  3.791915156825920e-09 \\
};
\addlegendentry{$v$, $e_{\infty}$}

\addplot [color=teal, dashed, line width=1.0pt, mark=diamond*, mark options={solid}, mark size=2.0pt]
table[row sep=crcr]{
7.230824819348436e-03  7.296986279041190e-05 \\
5.787060338230405e-03  2.383511565353599e-05 \\
3.877843685051529e-03  3.586406187773742e-06 \\
2.103067914659541e-03  1.599986141322898e-07 \\
9.087673506033234e-04  1.936495993938001e-09 \\
6.907543842015960e-04  5.216472121637508e-10 \\
};
\addlegendentry{$v$, $e_{2}$}

\addplot [color=black, dashed, line width=1.0pt]
table[row sep=crcr]{
7.230824819348436e-03  9.863859462878705e-04 \\
5.787060338230405e-03  3.238901408087152e-04 \\
3.877843685051529e-03  4.375817186421628e-05 \\
2.103067914659541e-03  2.052930107464759e-06 \\
9.087673506033234e-04  3.092938283772738e-08 \\
6.907543842015960e-04  7.847402454624339e-09 \\
};
\addlegendentry{$O(h^{5})$}

\end{axis}
\end{tikzpicture}
\end{subfigure}
\caption{Relative errors as functions of $h$ for computing \cref{eq:sysgreen3rd} using asymptotic expansions retaining terms through $\delta^{2}$. The system coefficients are $\alpha_2 = -1.31$, $\vb = [-7.23,4.3]$, $\alpha_3 =1$ and $\beta_3=-1$, with $\alpha_1=2$ for the left plot and $\alpha_1=200$ for the right plot.}
  \label{fig:5thorder}
\end{figure}

\section{Conclusions}
We have presented a framework for the evaluation of layer and volume potentials for strongly elliptic PDEs in complex geometries without requiring an explicit Green's function in physical space. The central idea is that an Ewald-type splitting of the Fourier symbol of the underlying differential operator decomposes the solution into a smooth, nonlocal component $\uh$, amenable to FFT-based evaluation, and a spatially localized component $\ul$, evaluated via explicit asymptotic expansions in powers of $\sqrt{\veps}$. The coefficients in these expansions are given in terms of local geometric quantities and derivatives of the source data at the closest boundary point, requiring no singular quadrature or kernel-specific constructions. The derivation proceeds entirely in the Fourier domain and requires only the symbol of the operator, making it directly applicable to coupled systems for which no closed-form Green's function exists.

For the Poisson equation in two and three dimensions, the resulting expansions recover the known structure of the Green's function-based formulas of \cite{lghtwght_ps2026}. This provides a check on the Fourier-domain derivation in the cases where the classical kernels are available. The extension to coupled strongly elliptic systems in \cref{sec:coupsys} is new, with the coupling between components entering naturally through the expansion of the matrix exponential. This is the main point of working with the symbol, because once it is known, the same localization and moment-reduction procedure applies even when no closed-form Green's function is available. The numerical experiments in \cref{sec:numerical_results} support this construction, showing fifth-order convergence when terms through $\delta^2$ are retained, including in the large-coupling test where the error level increases but the observed convergence rate is unchanged.

An important direction for future work is the integration of the present framework with the DMK algorithm \cite{dmk2024}, which would enable implementation on adaptive data structures with the splitting parameter $\veps$ chosen locally, and yield fully adaptive, high-order solvers for PDEs in complex geometries. Natural extensions include parabolic potential theory and data-driven settings where the impulse response is available through measurements rather than closed-form expressions. We note that when the free-space Green's function is known, the method \cite{beale2025} has advantages and a hybrid method is under investigation. Finally, splits other than those obtained from parabolic regularization will be investigated and reported in future work.

\section{Acknowledgments}
The author would like to thank Leslie Greengard (Courant Institute of Mathematical Sciences, New York University, and Center for Computational Mathematics, Flatiron Institute), Charles Epstein (Center for Computational Mathematics, Flatiron Institute), Shidong Jiang (Center for Computational Mathematics, Flatiron Institute), and Manas Rachh (Department of Mathematics, Indian Institute of Technology Bombay) for helpful discussions. The author gratefully acknowledges support from the Knut and Alice Wallenberg Foundation under grant 2020.0258.

\bibliographystyle{abbrv}
\bibliography{references}

\end{document}